\documentclass[12pt]{amsart}

\usepackage{a4wide,amssymb,bm,comment,mathbbol,mathrsfs,stmaryrd,times}

\usepackage[arrow,matrix,curve]{xy}\SilentMatrices
\newdir{ >}{{}*!/-1em/\dir{>}}
\def\xymat#1{\begin{aligned}\xymatrixcolsep{3pc}\xymatrix{#1}\end{aligned}}

\def\brk#1{\left\langle#1\right\rangle}
\def\set#1{\left\{#1\right\}}

\let\cl\brk
\def\acl#1{\hat{\mathstrut#1}}
\def\cref#1#2#3{(#2\mid#3)_{#1}}

\def\restr#1#2{\left.#1\right|_{#2}}
\def\1{^{-1}}
\def\2{{\scriptstyle{\frac12}}}
\def\com#1{\left[#1,#1\right]}

\let\lim\varprojlim

\let\eps\varepsilon
\let\Ph\varPhi
\let\le\leqslant
\let\ge\geqslant

\makeatletter

\let\cd@equation=\cd@subsubsection

\def\thmhead#1#2#3{%
 {\the\thm@notefont \thmnumber{
    \@upn{#2}}.}
   \thmname{#1}%
   \thmnote{ {\the\thm@notefont(#3)}}}

\makeatother

\def\theodef#1{\newtheorem{#1}[subsection]{#1}}

\theodef{Lemma} \theodef{Corollary} \theodef{Proposition} \theodef{Theorem}

\theoremstyle{definition}

\theodef{Definition}

\theodef{Remark} \theodef{Example} \theodef{Examples}

\DeclareMathAlphabet{\bsf}{OT1}{cmss}{bx}{n}

\def\xto#1{\xrightarrow[]{#1}}

\newcommand\Gr{\bsf{Groups}}
\newcommand\Ab{\bsf{Ab}}

\newcommand\Niq{{\bsf{Niq}}}
\newcommand\Spl{{\bsf{Spl}}}
\newcommand\Splab{{\Spl\ab}}
\newcommand\SpN{{\Spl(\Niq)}}
\newcommand\SpNN{{\Spl(\Niqssim)}}
\newcommand\Nil{{\bsf{Nil}}}

\newcommand\qu{\mathsf{Quad}}
\newcommand\wqu{\mathsf{wQuad}}

\newcommand\qw{{\mathsf Q}}

\newcommand\ab{^{\mathrm{ab}}}
\newcommand\op{^{\mathrm{op}}}
\newcommand\nilp{^{\mathrm{nil}}}

\def\Nilab{{\Nil\ab}}
\def\Niqab{{\Niq\ab}}
\def\Nilsim{{\Nil^\sim}}
\def\Niqsim{{\Niq^\sim}}
\def\Niqssim{{\Niq^\approx}}

\DeclareMathOperator\Id{\mathsf{Id}} \DeclareMathOperator\im{\mathsf{Im}}
\DeclareMathOperator\cok{\mathsf{Coker}} \DeclareMathOperator\Ker{\mathsf{Ker}}
\DeclareMathOperator\Ext{\mathsf{Ext}} \DeclareMathOperator\Hom{\mathsf{Hom}}
\DeclareMathOperator\centr{{\mathsf Z}} 
\DeclareMathOperator\sym{S}

\let\tp\otimes
\let\x\times
\let\ox\otimes

\let\into\rightarrowtail
\let\onto\twoheadrightarrow

\newcommand\pol{{\mathrm P}}

\newcommand\bC{{\boldsymbol C}}
\newcommand\bE{{\boldsymbol E}}
\newcommand\bZ{{\mathbf Z}}

\begin{document}

\title
{Quadratic envelope of the category of class two nilpotent groups}

\author{M. Jibladze}

\author{T. Pirashvili}

\address{Razmadze Mathematical Institute, M.~Alexidze st. 1, Tbilisi 0193,
Georgia}

\thanks{Research supported by the RTN Network HPRN-CT-2002-00287. The second
author supported by the Humboldt Foundation and the Deutsche
Forschungsgemeinschaft}

\maketitle

\section*{Introduction}

In this paper we introduce a somewhat surprising extension of the category of
class two nilpotent groups. It has the same objects but, unlike the latter,
its morphisms are closed under pointwise addition of maps. At the same time
the class of its morphisms is much smaller than the class of all maps between
groups. In fact, the morphisms are quadratic maps of very special kind which
we call q-maps. By definition a map $f:G\to H$ is called a q-map, if the
expression $(x\mid y)_f=-(f(x)+f(y))+f(x+y)$ lies in the commutator subgroup
of $H$ and is linear in $x$ and $y$. Any homomorphism is a q-map, but as we
said, the sum and composite of two q-maps is still a q-map and therefore one
obtains the category $\Niq$, with the objects all nilpotent groups of class
two and morphisms all q-maps between them. The advantage of $\Niq$ is the
fact that hom's in $\Niq$ are still nil$_2$-groups. Composition in $\Niq$ is
left distributive, but not right distributive. Actually $\Niq$ is an example
of a right quadratic category in the sense of \cite{BHP}. Since the category
$\Niq$ contains more morphisms than $\Nil$, two nonisomorphic groups might be
isomorphic as objects of $\Niq$. Thus the classification problems (say of
finite groups) in $\Niq$ are easier (but still highly nontrivial) than the
corresponding problems in $\Nil$.

We also indicate an approach to such classification questions using the
notion of linear extension of categories from \cite{BW}. Namely, we construct
several linear extensions connecting the category $\Niq$ to some simpler
categories, among them some additive ones which might be susceptible to
representation-theoretic classification methods. The point is that  a linear
extension induces bijection on isomorphism classes of objects.

Here is a short description of the contents of the paper. We begin by
reproducing in a maximally elementary way some known facts about the category
of class two nilpotent groups in section \ref{ni2}. In section \ref{quadra} we
present two variants of the notion of quadratic map for non-abelian groups and
investigate some basic properties of such maps. Then in section \ref{qmaps} we
introduce the new class of q-maps, lying strictly between homomorphisms and
quadratic maps and obtain various key properties of this class. Then in the
central section \ref{niq} using the q-maps we introduce the category $\Niq$
and give some of its features.

We then continue the study of $\Niq$ using linear extensions of categories.
In the next section \ref{nillinext} we recall this notion and exhibit the
category $\Nil$ of class two nilpotent groups and homomorphism as a linear
extension of a simpler category $\Nilsim$ that up to equivalence can be described in
terms of 2-cohomology classes of abelian groups. Then in section \ref{niqext}
we do a similar thing with $\Niq$ in place of $\Nil$; this time the simpler
category $\Niqsim$ is even additive, unlike $\Nilsim$, and moreover is itself a
linear extension of an even smaller additive category $\Niqssim$.

In the next section \ref{qsplit} we introduce a particular class of
nil$_2$-groups which we call q-split. This class seems to be a simplest
nontrivial one admitting classification modulo isomorphism in $\Niq$ in terms
of abelian groups. At the same time, it is quite rich, and smallest examples
of non-q-split groups are not quite trivial.

In section \ref{liecase} we exhibit an analog of the notion of q-map and the
category $\Niq$ for Lie algebras and prove that in the uniquely 2-divisible
situation the classical Maltsev correspondence between nil$_2$ groups and Lie
algebras extends to q-maps. This fact has some consequences for the
classification questions in view of further linear extensions on the Lie
algebra side. Finally in the last section \ref{obstr} we, using methods of
nonabelian cohomology, construct an obstruction to lifting homomorphisms to
q-maps and in particular find an obstruction for a nil$_2$-group to be
q-split.

\section{The category $\Nil$}\label{ni2}

The material in this section is well known and included for convenience of the
reader and to compare with what follows next. 


We fix some notation. Groups will be written additively.  For a group $G$
and elements $a,b\in G$ we let $[a,b]=-a-b+a+b$ be the commutator of $a$ and
$b$. If $G_1$ and $G_2$ are subgroups of $G$, then $[G_1,G_2]$ denotes the
subgroup generated by elements $[a,b]$, where $a\in G_1$ and $b\in G_2$. An
element $a\in G$ is called \emph{central} if $[a,x]=0$ for all $x\in G$. We
denote by $\centr(G)$ the center of $G$, which is the subgroup consisting of
all central elements of $G$.

For any group $G$ we denote by $G\ab$ the abelianization of $G$, that is,
the quotient
$$
G\ab:=G/[G,G].
$$
For an element $x\in G$ we let $\hat{x}$ denote the class of $x$ in $G\ab$. For
any abelian group $A$ one denotes by $\Lambda ^2(A)$ the \emph{second exterior
power} of $A$, which is the quotient of $A\tp A$ by the subgroup generated by
elements of the form $a\tp a$, $a\in A$.
 
A subgroup $A$ of a group $G$ is called \emph{central} if $[G,A]=0$, in other
words $A\subset\centr(G)$. A short exact sequence of groups 
$$
\bE=\left(0\to A\xto iG\xto pQ\to0\right)
$$
is called  a \emph{central extension} of $Q$ by $A$ if $i(A)$ is a central
subgroup of $G$. We refer to, e.~g., \cite{stammbach} for details on the
relationship between central extensions and the second cohomology.

A group $G$ is of \emph{nilpotence class two}, or is a \emph{nil$_2$-group},
if all triple commutators of $G$ vanish, $[[G,G],G]=0$, i.~e. one has
$[G,G]\subseteq\centr(G)$.

The smallest nonabelian groups of nilpotence class two are the quaternion
group $Q_8$ and the dihedral group $D_4=\bZ/4\bZ\ltimes\bZ/2\bZ$, both of
order 8. We denote by $\Nil$ the category of groups of nilpotence class two.

\begin{Lemma}\label{na+nb} For any $G\in\Nil$ one has:
\begin{itemize}
\item[i)] There is a well-defined homomorphism $\Lambda^2(G\ab)\to G$ given by
$\acl a\wedge\acl b\mapsto[a,b]$.

\item[ii)] For any $a,b\in G$ one has $[a,b]=a+b-a-b$.

\item[iii)] One has the inclusion $[G,G]\subset\centr(G)$.

\item[iv)] For any $a,b\in G$ and any $n\in\bZ$ one has
$$
na+nb=n(a+b)+\frac{n(n-1)}2[a,b].
$$
\end{itemize}
\end{Lemma}\qed

The inclusion functor $\Nil\subset\Gr$ has a left adjoint, given by
$$
G\mapsto G\nilp:=G/[[G,G],G].
$$
Since left adjoints preserve all existing colimits, one can obtain coproducts
in $\Nil$ as $(-)\nilp$ of coproducts in $\Gr$. But in fact coproducts in
$\Nil$ are much easier to construct directly than those in $\Gr$. Namely one
has


\begin{Proposition}\label{nilcopr}
For two nil$_2$-groups $G$, $H$ let $G\vee H$ be the set $G\ab\!\ox\!H\ab\x G\x
H$. The equalities
\begin{align*}
&(\xi,g,h)+(\xi',g',h')=(\xi+\xi'-\acl{g'}\ox\acl h,g+g',h+h'),\\
&-(\xi,g,h)=(-\xi-\acl g\ox\acl h,-g,-h),\\
&0=(0,0,0)
\end{align*}
equip this set with a nil$_2$-group structure such that there is a central
extension
$$
0\to G\ab\ox H\ab\to G\vee H\to G\x H\to 0.
$$
Moreover the maps $i_G:G\to G\vee H$, $i_H:H\to G\vee H$ given by
$i_G(g)=(0,g,0)$ and $i_H(h)=(0,0,h)$ form a coproduct diagram in $\Nil$.
\end{Proposition}

\begin{proof}
The map $(G\x H)\x(G\x H)\to G\ab\ox H\ab$ given by $((g,h),(g',h'))\mapsto
-\acl{g'}\ox\acl h$ is easily seen to be a 2-cocycle, so it indeed defines a central
extension as above, and the indicated maps are clearly homomorphisms. One
calculates
$$
[(\xi,g,h),(\xi',g',h')]
=(\hat{g}\ox\acl{h'}-\acl{g'}\ox\hat{h},[g,g'],[h,h']);
$$
in particular, it follows that the elements $(0,[g,g'],0)$ and $(0,0,[h,h'])$,
along with $(\xi,0,0)$, are central in $G\vee H$, so that the latter is a
nil$_2$-group.

Moreover in $G\vee H$ one obviously has the identities
$$
(\xi,g,h)=(\xi,0,0)+(0,g,0)+(0,0,h)
$$
and
$$
(\hat{g}\ox\acl h,0,0)=[(0,g,0),(0,0,h)].
$$
Hence if we want to extend some homomorphisms $u:G\to X$, $v:H\to X$ with
$X\in\Nil$ to a homomorphism $(u,v):G\vee H\to X$ along $i_G$ and $i_H$, by the
above identities we have a unique choice, namely to put
$$
(\xi,g,h)\mapsto[u,v](\xi)+u(g)+v(h),
$$
where $[u,v]:G\ab\ox H\ab\to X$ is determined by $[u,v](\hat{ g}\ox\acl
h)=[u(g),v(h)]$. Since $X$ is in $\Nil$, the expression $[u(g),v(h)]$ factors
through $G\ab\x H\ab$ and is bilinear, so we indeed have a correctly defined
map $G\vee H\to X$. Then using the fact that the elements $[u(g),v(h)]$ are
also central in $X$, it is easy to see that this map is in fact a homomorphism.
\end{proof}


The forgetful functor $\Nil\to\bsf{Sets}$ has a left adjoint, whose value on a
set $S$ is known as \emph{the free nilpotent group of class two generated by}
$S$ and is denoted by $\bZ_\Nil[S]$. If $F_S$ is the free group spanned by $S$,
then $\bZ_\Nil[S]=(F_S)\nilp$. Moreover since left adjoints preserve
coproducts, and $S$ is the coproduct of $S$ copies of a singleton in
$\bsf{Sets}$, one has
$$
\bZ_\Nil[S]=\bigvee_S\bZ
$$
in $\Nil$. Using Proposition \ref{nilcopr}, we obtain the following  particular case of the
famous result of Witt, which asserts that the graded Lie ring obtained by the
lower central series of a free group is a free Lie ring.

\begin{Corollary}\label{wittori}
For a free nil$_2$-group $G$ one has the following central extension
$$
0\to \Lambda^2(G\ab)\to G\to G\ab\to0
$$
\end{Corollary}

\begin{proof}
It suffices to prove the lemma for $G=\bZ_\Nil[S]$ with $S$ finite. Indeed
every $S$ is a directed colimit of its finite subsets, all functors under
consideration preserve colimits, and a directed colimit of short exact
sequences is short exact.

For finite $S$ we use induction on the number of elements, the case of one
element, i.~e. $G=\bZ$, being trivially true. In other words, we have to show
that if the above sequence is short exact for $G$ and $H$, then it also is for
$G\vee H$. Now this is clear from the following diagram with exact rows and
columns
$$
\xymat{
       &0\ar[d]                                     &0\ar[d]\\
       &G\ab\oplus H\ab\ar[d]\ar[r]^=               &G\ab\oplus H\ab\ar[d]\\
       &\Lambda^2((G\vee H)\ab)\ar[r]\ar[d]         &G\vee H\ar[r]\ar[d]&(G\vee H)\ab\ar[d]_\cong\ar[r]&0\\
0\ar[r]&\Lambda^2(G\ab)\x\Lambda^2(H\ab)\ar[r]\ar[d]&G\x H\ar[r]\ar[d]  &G\ab\x H\ab\ar[r]&0\\
       &0                                           &0
}
$$
taking into account that for any $G$, $H$ the canonical homomorphism
$G\ab\oplus H\ab\to(G\vee H)\ab$ is an isomorphism since $(-)\ab$ is a left
adjoint and thus preserves all colimits.
\end{proof}

\section{Quadratic maps between nonabelian groups}\label{quadra}

Let $G$ and $H$ be arbitrary groups. Call a map
$f:G\to H$ \emph{weakly quadratic} if for any $a,b\in G$ the
\emph{cross-effect}
$$
\cref fab:=-(f(a)+f(b))+f(a+b)
$$
commutes with $f(c)$ for all $c\in G$ and is linear in $a$ and $b$. Thus we
have
$$
f(a+b)=f(a)+f(b)+\cref fab,
$$
and the equalities
$$
\begin{array}{rl}
&\cref f{a_1+a_2}b=\cref f{a_1}b+\cref f{a_2}b,\\
&\cref fa{b_1+b_2}=\cref fa{b_1}+\cref fa{b_2},\\
&\cref fab+f(c)=f(c)+\cref fab
\end{array}
$$
hold for any $a,a_1,a_2,b_1,b_2,b,c\in G$.

A weakly quadratic map $f:G\to H$ is \emph{quadratic} \cite{BHP} if in fact
$\cref fab\in\centr(H)$ for all $a,b\in G$.

Obviously every weakly quadratic map to an abelian group is quadratic. We
denote the set of all weakly quadratic maps from $G$ to $H$ by $\wqu(G,H)$ and
that of quadratic maps by $\qu(G,H)$. It is clear that a map $f:G\to H$ is a
homomorphism iff $\cref f--=0$. Thus
$$
\Hom(G,H)\subseteq\qu(G,H)\subseteq\wqu(G,H).
$$

\begin{Lemma}\label{qwadasxvistvisebebi}
For $f\in\wqu(G,H)$ the following assertions are true:
\begin{itemize}
\item[i)] The cross-effect yields a well-defined homomorphism $\cref
f--:G\ab\ox G\ab\to H$.

\item[ii)] $f(0)=0$.

\item[iii)] $f(-a)=-f(a)+\cref faa$.

\item[iv)] If $c\in[G,G]$, then for any $a\in G$ one has $f(a+c)=f(a)+f(c)$. In
particular the restriction of $f$ to the commutator subgroup is a homomorphism.

\item[v)] For any $a,b\in G$ one has
$$
f([a,b])=-f(b+a)+f(a+b)=[f(a),f(b)]+\cref fab-\cref fba.
$$

\item[vi)] For any $a,b,c\in G$ one has $f([a,[b,c]])=[f(a),[f(b),f(c)]]$.
\end{itemize}
\end{Lemma}

\begin{proof}
i) Since the elements $\cref fab$ commute with everything in the image of $f$,
they centralize the subgroup generated by this image. But they belong to this
subgroup themselves, so commute with each other. Thus for each $a\in G$ the map
$\cref fa-:G\to H$ is a homomorphism with abelian image, hence it factors
trough $G\ab$. Similarly for $\cref f-b$.

ii) $0=\cref f00=-f(0)-f(0)+f(0)=f(0)$.

iii) By ii),
$$
0=f(0)=f(a-a)=f(a)+f(-a)+\cref fa{-a}
$$
and the statement follows.

iv) By i), $\cref fac=0$.

v) We have
\begin{flalign*}
f([a,b])&=f(-(b+a)+a+b)\\
&=f(-(b+a))+f(a+b)+\cref f{-(b+a)}{a+b}\\
&=-f(b+a)+\cref f{b+a}{b+a}+f(a+b)+\cref f{-(b+a)}{a+b}&\text{(by iii))}\\
&=-f(b+a)+f(a+b)&\text{(by i))}\\
&=-(f(b)+f(a)+\cref fba)+f(a)+f(b)+\cref fab\\
&=[f(a),f(b)]+\cref fab-\cref fba.
\end{flalign*}

vi) We have
\begin{flalign*}
f([a,[b,c]])&=[f(a),f([b,c])]+\cref fa{[b,c]}-\cref f{[b,c]}a&\text{(by v))}\\
&=[f(a),f([b,c])]&\text{(by i))}\\
&=[f(a),[f(b),f(c)]+\cref fbc-\cref fcb]&\text{(by v))}\\
&=[f(a),[f(b),f(c)]]&\text{(by i))}.
\end{flalign*}
\end{proof}

\begin{Corollary}
Let $f:G\to H$ be a weakly quadratic map. If $H$ is a nilpotent group of class
two, then $f$ factors through $G\nilp=G/[G,[G,G]]$. Thus
\begin{align*}
\wqu(G,H)&\cong\wqu(G\nilp,H)\\
\qu(G,H)&\cong\qu(G\nilp,H).
\end{align*}
\end{Corollary}
\begin{proof}
Indeed, if $c\in[G,[G,G]]$ then $f(c)=0$ thanks to vi) of Lemma
\ref{qwadasxvistvisebebi}. Thus $f(a+c)=f(a)$ by iv) of Lemma
\ref{qwadasxvistvisebebi}.
\end{proof}

The set of quadratic maps for nilpotent groups of class two has some remarkable
properties. First of all unlike $\Hom(G,H)$ or $\wqu(G,H)$ the set $\qu(G,H)$
is a group with respect to the  pointwise addition of maps. This is the subject
of the following Lemma.

\begin{Lemma}\label{jamiscross}
Let $G$ be a group and let $H$ be a nilpotent group of class two. If the maps
$f,g:G\to H$ are quadratic, then  $f+g$ and $-f$ are also quadratic and
\begin{align*}
\cref{f+g}ab&=\cref fab+\cref gab+[f(b),g(a)],\\
\cref{-f}ab&=[f(b),f(a)]-\cref fab.
\end{align*}
\end{Lemma}

\begin{proof}
The above formul\ae\ for $\cref{f+g}--$ and $\cref{-f}--$ can be easily
checked. Since the commutators are central, it remains to show that
$[f(b),g(a)]$ is linear in $a$ and $b$ for any quadratic $f$ and $g$. But this
is clear, because $[-,-]$ is central, bilinear and vanishes on central
elements.
\end{proof}

\begin{Example}[\cite{square}]\label{uniqu}
For any group $G$ there exists a \emph{universal weakly quadratic map}
$p_2:G\to\pol_2G$ such that for any other weakly quadratic map $q:G\to H$ there
is a unique homomorphism $f_q:\pol_2G\to H$ with $q=f_qp_2$. Thus
$$
\wqu(G,H)\cong\Hom(\pol_2G,H).
$$
One defines $\pol_2G$ by the pullback square
$$
\xymatrix{
\pol_2G\ar@{->>}[r]\ar@{ >->}[d]\ar@{}[dr]|<{\Bigg\lrcorner}&G\ar@{ >->}[d]^{\mathrm{diagonal}}\\
G\vee G\ar@{->>}[r]&G\x G. }
$$
Thus by Proposition \ref{nilcopr} there is a central extension
\begin{equation}\label{314Io}
0\to G\ab\ox G\ab\xto\iota\pol_2G\xto\pi G\to 0.
\end{equation}
and $\pol_2G$ is isomorphic to the set $G\ab\!\ox\!G\ab\x G$ with the group
structure given by
$$
(\xi,g)+(\xi',g')=(\xi+\xi'-\hat{g}\ox\hat{g'},g+g').
$$

The universal weakly quadratic map $p_2:G\to\pol_2G$ is given by
$p_2(g)=(0,g)$. Indeed in $\pol_2G$ one then has
$$
(\hat{x}\ox\hat{y},0)=-((0,x)+(0,y))+(0,x+y)=\cref{p_2}xy
$$
and
$$
(\xi,g)=(\xi,0)+p_2(g),
$$
so to factor a weakly quadratic map $q:G\to H$ through $p_2$ via a homomorphism
$f_q:\pol_2G\to H$ one is forced to put
$$
f_q\left(\xi,g\right)=\cref q--(\xi)+q(g).
$$
One then checks easily that this indeed gives the required factorization.
 
Note that the universal weakly quadratic map $p_2:G\to\pol_2G$ is not only
weakly quadratic but actually also quadratic.
\end{Example}
\begin{Lemma} For any $G\in \Niq$ and $A\in \Ab$ one has an exact sequence:
$$0\to \Hom(G,A)\to \qu(G,A)\to \Hom(G\ab\tp G\ab,A)\to H^2(G,A)$$
where the last homomorphism is given by $f\mapsto f_*([G])$. Here
$f:G\ab\tp G\ab\to A$ is a homomorphism and $[G]\in H^2(G,G\ab\tp G\ab)$ 
is the class represented by the central extension \eqref{314Io}.
\end{Lemma}

\begin{proof}
The result is an immediate consequence of the 5-term exact sequence in group
cohomology (see for example  \cite{stammbach}, page 15) applied to the
central extension  \eqref{314Io} and from the fact that for abelian $A$ one has
$\Hom(\pol_2G,A)\cong \qu(G,A)$.
\end{proof}

\begin{Lemma}\label{prods}
For any groups $(G_i)_{i\in I}$ and $H$ one has natural bijections
\begin{align*}
\qu(H,\prod_iG_i)\approx\prod_i\qu(H,G_i).
\end{align*}
If moreover $H\in\Nil$ then there is a central extension
$$
0\to\Hom(G_1\ab\ox G_2\ab,\centr(H)) \xto\alpha\qu(G_1\x
G_2,H)\to\qu(G_1,H)\x\qu(G_2,H)\to0
$$
where $(\alpha(\xi))(g_1,g_2)=\xi(\hat g_1,\hat g_2)$ for
$\xi\in\Hom(G_1\ab\ox G_2\ab,\centr(H))$ and $g_k\in G_k$, $k=1,2$.
\end{Lemma}

\begin{proof}
The first assertion is clear. For the second, take any elements
$f_k\in\qu(G_k,H)$, $k=1,2$. Then the composite maps $f_kp_k$ are again
quadratic, where $p_k:G_1\x G_2\to G_k$ are the projections. Thus
$f=f_1p_1+f_2p_2:G_1\x G_2\to H$ is a quadratic map. It is clear that
$fi_k=f_k$, where $i_k:G_k\to G_1\x G_2$ are the standard inclusions. This
shows that the map $\qu(G_1\x G_2,H)\to\qu(G_1,H)\x\qu(G_2,H)$ is surjective.
Let us compute the kernel of the latter homomorphism. Take an $f$ from the
kernel. Then $f:G_1\x G_2\to H$ is a quadratic map such that
$f(g_1,0)=0=f(0,g_2)$ for all $g_k\in G_k$. Define $\xi:G_1\ab\x
G_2\ab\to\centr(H)$ by $\xi(\acl{g_1},\acl{g_2}):=\cref f{(g_1,0)}{(0,g_2)}$.
Then one has
$$
f(g_1,g_2)=f((g_1,0)+(0,g_2))=\cref
f{(g_1,0)}{(0,g_2)}=\xi(\hat{g_1},\hat{g_2})
$$
and the lemma follows.
\end{proof}

In the rest of the paper we will assume that all groups under consideration are
nilpotent of class two.

\begin{Lemma}\label{komposiciiscros}
Let $f:G\to H$ be a weakly quadratic map. For any homomorphism $h:G_1\to G$ the
composite $fh:G_1\to H$ is also weakly quadratic and
$$
\cref{fh}ab=\cref f{h(a)}{h(b)}, \ \ a,b\in G_1;
$$
moreover if $f$ is quadratic then so is $fh$.

For any homomorphism $g:H\to H_1$ the composite $gf:G\to H_1$ is also weakly
quadratic and
$$
\cref{gf}ab=g(\cref fab).
$$
If moreover $f$ is quadratic then $gf$ will be quadratic provided $g$ carries
central elements to central elements.
\end{Lemma}\qed

Thus for any $N\in\Nil$, one obtains functors
\begin{align*}
\wqu(-,N)&:\Nil\op\to\bsf{Sets},\\
\qu(-,N)&:\Nil\op\to\Nil,\\
\wqu(N,-)&:\Nil\to\bsf{Sets}.
\end{align*}
In fact by Example \ref{uniqu} the last functor is representable, i.~e. one has
$$
\wqu(N,-)\approx\Hom_\Nil(\pol_2N,-).
$$
However the mapping $\qu(N,-)$ is NOT functorial.

\begin{Examples}\label{qwasxvebi}

i) For a fixed $n\in\bZ$, consider the map $n:G\to G$ given by $a\mapsto na$.
Then
$$
\cref nab=-\frac{n(n-1)}{2}[a,b].
$$
Thus $n\in\qu(G,G)$.

ii) Let $+:G\x G\to G$ be the map given by $(a,b)\mapsto a+b$. Then
$$
\cref+{(a,b)}{(c,d)}=[c,b].
$$
In particular $+\in\qu(G\x G,G)$.

iii) For any elements $a\in G$ and $b\in\centr(G)$ we put
$$
f_{a,b}(n)=na+\frac{n(n-1)}2b.
$$
The map $f_{a,b}:\bZ\to G$ is a quadratic map with $\cref{f_{a,b}}nm=nmb$ for
any $n,m\in\bZ$. We claim that any quadratic map $f:\bZ\to G$ is of this form.
Indeed, one puts $a=f(1)$, $b=\cref f11$ and considers $g=f-f_{a,b}$. Then one
has $g(1)=0=\cref g11$. Since $\cref g--$ is bilinear it follows that $\cref
gnm=nm\cref g11=0$. Hence $g$ is a homomorphism and the condition $g(1)=0$
shows that $g=0$ and the claim is proved. One easily computes that
$$
f_{a,b}+f_{a',b'}=f_{a+a',b+b'+[a,a']}.
$$
Thus pointwise sum of quadratic maps $\bZ\to G$ is quadratic, so that
$\qu(\bZ,G)$ has a group structure and one has the following central extension:
$$
0\to\centr(G)\to\qu(\bZ,G)\xto{\mathsf{ev}(1)}G\to0
$$
where $\mathsf{ev}(1)(f)=f(1)$. A 2-cocycle $G\x G\to\centr(G)$ corresponding
to this central extension is given by the commutator map.
\end{Examples}

Let us next investigate quadratic maps of the form $f:G_1\vee G_2\to H$. For
such a map, denote
$$
f_i=\restr f{G_i}:G_i\to H,
$$
$i=1,2$ and
$$
f_\ox=\restr f{G_1\ab\ox G_2\ab}:G_1\ab\ox G_2\ab\to H,
$$
where the inclusion $G_1\ab\ox G_2\ab\subset G_1\vee G_2$ is as in Proposition
\ref{nilcopr}. Since $G_1\ab\ox G_2\ab$ is contained in the commutator subgroup
of $G_1\ox G_2$, the map $f_\ox$ is a homomorphism, and its image lies in the
center of $H$ (by v) of Lemma \ref{qwadasxvistvisebebi}). As for $f_i$, they are
quadratic maps. Since every element of $G\vee H$ has the form
$(\xi,g_1,g_2)=(\xi,0,0)+(0,g_1,0)+(0,0,g_2)$ with $\xi\in G_1\ab\ox G_2\ab$
and
$$
f(\xi,g_1,g_2)=f_\ox(\xi)+f_1(g_1)+f_2(g_2) +\cref f{\hat{g_1}}{\hat{g_2}},
$$
it follows that $f$ is uniquely reconstructed from the maps $f_\ox$, $f_1$,
$f_2$ and the homomorphism
$$
\restr{\cref f--}{G_1\ab\ox G_2\ab}:G_1\ab\ox G_2\ab\to\centr(H),
$$
which we will denote by $\acl f$.

Conversely, for any given maps
$$
f_i\in\qu(G_i,H),\ \ i=1,2,\ \ \ f_\ox,\acl f\in\Hom(G_1\ab\ox
G_2\ab,\centr(H))
$$
define the map $f:G_1\vee G_2\to H$ by
$$
f(\xi,g_1,g_2)=f_\ox(\xi)+f_1(g_1)+f_2(g_2) +\hat{f}(\hat{g_1}\ox\hat{g_2}).
$$
Then
\begin{align*}
f(&(\xi,x_1,x_2)+(\eta,y_1,y_2))=f(\xi+\eta-\hat{y_1}\ox\hat{x_2},x_1+y_1,x_2+y_2)\\
=f_\ox(&\xi+\eta-\hat{y_1}\ox\hat{x_2})+f_1(x_1+y_1)+f_2(x_2+y_2)
+\acl f((\hat{x_1}+\hat{y_1})\ox(\hat{x_2}+\hat{y_2}))\\
=f_\ox(&\xi)+f_\ox(\eta)-f_\ox(\hat{y_1}\ox\hat{x_2})\\
&+f_1(x_1)+f_1(y_1)+\cref{f_1}{x_1}{y_1}
+f_2(x_2)+f_2(y_2)+\cref{f_2}{x_2}{y_2}\\
&+\hat{f}(\hat{x_1}\ox\hat{x_2}) +\hat{f}(\hat{x_1}\ox\hat{y_2}) +\hat{
f}(\hat{y_1}\ox\hat{x_2})
+\hat{f}(\hat{y_1}\ox\hat{y_2})\\
=f(&\xi,x_1,x_2)+f(\eta,y_1,y_2)\\
&-f_\ox(\hat{y_1}\ox\hat{x_2})+[f_1(y_1),f_2(x_2)] +\cref{f_1}{x_1}{y_1}
+\cref{f_2}{x_2}{y_2} +\hat{f}(\acl{x_1}\ox\acl{y_2}) +\hat{
f}(\hat{y_1}\ox\hat{x_2}).
\end{align*}
It follows that $f$ is a quadratic map, so that indeed any choice of $f_\ox$,
$f_1$, $f_2$ and $\acl f$ as above is valid.

Now suppose given two quadratic maps $f,f':G_1\vee G_2\to H$. Then for their
sum clearly one has $(f+f')_i=f_i+f'\mathstrut_i$, $i=1,2$, and
$(f+f')_\ox=f_\ox+f'\mathstrut_\ox$. Moreover one calculates
\begin{multline*}
\widehat{f+f'}(\hat{g_1}\ox\hat{g_2})
=\cref{f+f'}{(0,g_1,0)}{(0,0,g_2)}\\
=\cref f{(0,g_1,0)}{(0,0,g_2)}+\cref{f'}{(0,g_1,0)}{(0,0,g_2)}+[f(0,0,g_2),f'(0,g_1,0)]\\
=\hat{f}(\hat{g_1}\ox\acl{g_2})+\acl{f'}(\hat{g_1}\ox\hat{g_2})+[f_2(g_2),f'\mathstrut_1(g_1)].
\end{multline*}
Thus identifying $f$ with the quadruple $(f_1,f_2,f_\ox,\hat{f})$ as above one
has
$$
(f_1,f_2,f_\ox,\acl
f)+(f'\mathstrut_1,f'\mathstrut_2,f'\mathstrut_\ox,\acl{f'})
=(f_1+f'\mathstrut_1,f_2+f'\mathstrut_2,f_\ox+f'\mathstrut_\ox,\acl
f+\acl{f'}-[f'\mathstrut_1,f_2]),
$$
where
$$
[f'\mathstrut_1,f_2](\hat{g_1}\ox\hat{g_2})=[f'\mathstrut_1(g_1),f_2(g_2)].
$$

We thus have proved
\begin{Lemma}\label{quac}
For any nil$_2$-groups $G_1,G_2,H$ there is a central extension
\begin{align*}
0\to\Hom(G_1\ab\ox G_2\ab,\centr(H))\to\qu(G_1\vee G_2,H)&\to
\Hom(G_1\ab\ox G_2\ab,\centr(H))\\
&\x\qu(G_1,H)\x\qu(G_2,H)\to0.
\end{align*}
A cocycle defining this extension is given by
\begin{align*}
((f_\ox,f_1,f_2),(f'\mathstrut_\ox,f'\mathstrut_1,f'\mathstrut_2))\mapsto
\alpha((f_\ox,f_1,f_2),(f'\mathstrut_\ox,f'\mathstrut_1,f'\mathstrut_2))
:G_1\ab&\ox G_2\ab\to\centr(H),\\
\hat{g_1}&\ox\hat{g_2}\mapsto[f_2(g_2),f'\mathstrut_1(g_1)].
\end{align*}
\end{Lemma}\qed
\begin{Corollary} Let $G$ be a free nil$_2$-group on $x_1,...,x_n$. Then for
any nil$_2$-group $H$ and any elements $a_1,...,a_n\in H$, 
$a_{ij}, b_{ij}\in \centr(H)$, $i<j$ there exists a unique
quadratic map $f:G\to H$ such that
$$
\begin{array}{rll}
f(x_i)&=a_i,&1\le i\le n,\\
f([x_i,x_j])&=a_{ij},&i<j,\\
\cref f{x_i}{x_j} &= b_{ij},&i<j.
\end{array}
$$
\end{Corollary}\qed

\section{$\mathrm q$-maps}\label{qmaps}

The last identity of Lemma \ref{na+nb} suggests the following
definition:
\begin{Definition}\label{q}
A weakly quadratic map $f:G\to H$ between nil$_2$-groups is a \emph{q-map} if
one has $\cref fab\in[H,H]$ for all $a,b\in G$.
\end{Definition}
 We denote by $\qw(G,H)$ the collection
of all q-maps from $G$ to $H$, so that
$$
\Hom(G,H)\subseteq\qw(G,H)\subseteq\qu(G,H).
$$

\begin{Lemma}
The set $\qw(G,H)$ is a normal subgroup of $\qu(G,H)$. In particular  any
linear combination of homomorphisms is a q-map.
\end{Lemma}

\begin{proof} The first identity of Lemma  \ref{jamiscross} shows that
$\qw(G,H)$ is a subgroup of $\qu(G,H)$. By the same Lemma for any $f\in
\qu(G,H)$ and $g\in \qw(G,H)$ we have
$$\cref {f+g-f} ab=\cref g ab+ [fb,ga]-[gb,fa]$$
and the result follows.
\end{proof}

\begin{Lemma}\label{2x2}
A weakly quadratic map is in $\qw(G,H)$ iff its
composite with $H\onto H\ab$ is a homomorphism. In particular any q-map $f:G\to H$ yields
a well-defined homomorphism $f\ab:G\ab\to H\ab$ such that the diagram
$$
\xymatrix{
G\ar@{->>}[d]\ar[r]^f&H\ar@{->>}[d]\\
G\ab\ar[r]^{f\ab}&H\ab }
$$
commutes. 
\end{Lemma}

\begin{proof}
Indeed $\cref fab\in[H,H] $ iff the image of $\cref fab$ vanishes in $H\ab$.
\end{proof}

Obviously one has an embedding
$$
\qu(G,[H,H])\subset\qw(G,H)
$$
as a central subgroup.

\begin{Lemma}\label{qabelsi}
For an abelian group $H$ one has
$$
\qw(G,H)=\Hom(G,H)
$$
for any $G\in\Nil$.
\end{Lemma}

\begin{proof}
Since $[H,H]=0$, a map $f:G\to H$ is a q-map iff $\cref f--=0$.
\end{proof}

\begin{Examples}
The first two quadratic maps considered in Examples \ref{qwasxvebi} are
actually q-maps. Also the map
$$
\delta=i_1+i_2:G\to G\vee G
$$
is a q-map, with $\cref\delta xy=[i_1(y),i_2(x)]$.

On the other hand, the quadratic map $f_{a,b}:\bZ\to G$ associated to elements
$a\in G$ and $b\in\centr(G)$ as in iii) of Examples \ref{qwasxvebi} is a q-map iff
$b\in[G,G]$. Thus for any $G\in\Nil$ one has the following central extension:
$$
0\to[G,G]\to\qw(\bZ,G)\xto{\mathsf{ev}(1)}G\to0.
$$
A 2-cocycle $G\x G\to[G,G]$ corresponding to this central extension is given by
the commutator map.
\end{Examples}

Exactly as for  Lemma \ref{prods} one has

\begin{Lemma}\label{compcref}
For any groups $(G_i)_{i\in I}$, $H$ one has natural bijections
$$
\qw(H,\prod_iG_i)\cong\prod_i\qw(H,G_i).
$$
If moreover $H\in\Nil$ then there is a central extension
$$
0\to\qw(G_1\ab\ox G_2\ab,[H,H])\to\qw(G_1\x
G_2,H)\to\qw(G_1,H)\x\qw(G_2,H)\to0.
$$
\end{Lemma}\qed

Moreover one has exactly as in Lemma \ref{quac}

\begin{Lemma}
For any nil$_2$-groups $G_1,G_2,H$ there is a central extension
\begin{align*}
0\to\Hom(G_1\ab\ox G_2\ab,[H,H])\to\qw(G_1\vee G_2,H)&\to
\Hom(G_1\ab\ox G_2\ab,[H,H])\\
&\x\qw(G_1,H)\x\qw(G_2,H)\to0.
\end{align*}
A cocycle defining this extension is given by
\begin{align*}
((f_\ox,f_1,f_2),(f'\mathstrut_\ox,f'\mathstrut_1,f'\mathstrut_2))\mapsto
\alpha((f_\ox,f_1,f_2),(f'\mathstrut_\ox,f'\mathstrut_1,f'\mathstrut_2))
:G_1\ab&\ox G_2\ab\to[H,H],\\
\hat{g_1}&\ox\hat{g_2}\mapsto[f_2(g_2),f'\mathstrut_1(g_1)].
\end{align*}
In particular, if $G$ is a free nil$_2$-group on $x_1,\cdots,x_n$ then for any
nil$_2$-group $H$ and any elements $a_1,\cdots, a_n\in H$, $a_{ij}, b_{ij}\in [H,H], i<j$ there exists the unique q-map $f:G\to H$ such
that
$$
\begin{array}{rll}
f(x_i)&=a_i, &1\le i\le n,\\
f([x_i,x_j])&=a_{ij}, &i<j,\\
\cref f {x_i}{x_j} &= b_{ij}, &i<j.
\end{array}
$$
\end{Lemma}\qed

By Lemma \ref{2x2} any q-map $f:G\to H$ yields a homomorphism $f\ab:G\ab\to H\ab$. We now
associate two more  homomorphisms to any q-map.

\begin{Proposition} \label{2x3} Let $f:G\to H$ be a q-map. Then $f([G,G])\subset [H,H]$ and the
restriction of $f$ to $[G,G]$ yields a homomorphism $[f,f]:[G,G]\to [H,H]$, which fits in the
following commutative diagram
$$
\xymatrix{
0\ar[r]& [G,G]\ar[r]\ar[d]^{\com f}& G\ar[r]\ar[d]^f& G\ab\ar[r] \ar[d] ^{f\ab}&0\\
0\ar[r]& [H,H]\ar[r] & H\ar[r]& H\ab\ar[r]&0}
$$ 
Moreover, there exists a unique homomorphism 
$$
\beta(f):\Ker(f\ab)\to \cok([f,f])
$$
such that
$$
\beta(f)(\hat{a})=f(a)\mod f([G,G])
$$
for any $a\in f\1([H,H])$. 
Furthermore if $f$ is injective then $[f,f]$ and $\beta(f)$ are 
monomorphisms and if $f$ is surjective then $\beta(f)$ and $f\ab$ are epimorpisms. 
\end{Proposition}

\begin{proof}
If $f$ is a q-map then it follows from v) of Lemma \ref{qwadasxvistvisebebi} 
that $f[G,G]\subseteq[H,H]$. Hence by iv) of  Lemma \ref{qwadasxvistvisebebi}, 
$[f,f]:[G,G]\to [H,H]$ is a homomorphism and obviosly the diagram commutes. 
We claim that $\beta(f)$ is well-defined. One observes that if $\hat a_1=\hat a$
then $a_1=a+b$ with $b\in [G,G]$. It follows by iv) of  Lemma \ref{qwadasxvistvisebebi}
that $f(a_1)=f(a) \mod f([G,G])$ and the claim is proved. The rest is
just diagram chase.
\end{proof}

%
%

\section{The category $\Niq$}\label{niq}

In this section, our main character enters. This is the category
$\Niq$. The definition is based on the following result.

\begin{Proposition}\label{niqprops}
Any composite of q-maps is a q-map. More precisely, for q-maps $f:G\to H$ and
$g:G_1\to G$ the cross-effect of their composite is given by
$$
\cref{fg}ab=f(\cref gab)+\cref f{g(a)}{g(b)}, \ \ a,b\in G_1.
$$
\end{Proposition}

\begin{proof}
One has
\begin{flalign*}
fg(a+b)&=f(g(a)+g(b)+\cref gab)\\
&=f(g(a)+g(b))+f(\cref gab)&\text{(by iv) of Lemma \ref{qwadasxvistvisebebi})}\\
&=f(g(a))+f(g(b))+\cref f{g(a)}{g(b)}+f(\cref gab),
\end{flalign*}
which proves the equality above.
\end{proof}

Hence there is a well-defined category $\Niq$ whose objects are nil$_2$-groups
and morphisms are all q-maps between them. The hom-sets
$$
\Hom_{\Niq}(G,H)=\qw(G,H)
$$
are equipped with structures of nilpotent groups of class two. $\Nil$ is a
subcategory of $\Niq$, with the same objects. The hom-functor of $\Niq$ (with
values in sets) gives rise to a well-defined bifunctor
$$
\qw(-,-):\Nil\op\x\Nil\to\Nil.
$$
Moreover there are well-defined  functors $\Niq\to\Ab$ given respectively by
$G\mapsto G\ab$ and $G\mapsto [G,G]$.

Composition in $\Niq$ is left distributive,
$$
(f+f')g=fg+f'g,
$$
but not right distributive; rather it is right quadratic, in the following
sense. First  of all  one has
\begin{equation}
f(g+g')=fg+fg'+\cref fg{g'}, \ \ f\in \qw(G,H),\ \ g,g'\in\qw(G_1,G),
\end{equation}
where $\cref fg{g'}:G_1\to H$ is given by
\begin{equation}\label{komcr}
\cref fg{g'} (x)=  \cref f{g(x)}{g'(x)}.
\end{equation}
Secondly  $\cref fg{g'}$ lies in the center of $\qw(G_1,H)$ and it is bilinear
in $g$, $g'$ and quadratic in $f$
--- more precisely, one has
$$
\cref{f+f'}g{g'}(x)=\cref fg{g'}(x)+\cref{f'}g{g'}(x)+[fg'(x),f'g(x)].
$$

The category $\Niq$ possesses all products and both the inclusion
$\Nil\hookrightarrow\Niq$ and the forgetful functor $\Niq\to\bsf{Sets}$ respect
products.

Every object in $\Niq$ has a canonical internal group structure. However a
morphism in $\Niq$ is compatible with the corresponding internal group
structures iff it lies in $\Nil$, i.~e. is a homomorphism.

If $p_k:G_1\x G_2\to G_k$ are the standard projections and $i_k:G_k\to G_1\x
G_2$ are the standard inclusions then one has $p_ki_k=\Id_{G_k}$,
$i_1p_1+i_2p_2=\Id_{G_1\x G_2}$ and $p_2i_1=0$, $p_1i_2=0$. Therefore $\Niq$ is
a right quadratic category in the terminology of \cite{BHP}. Trivial groups are
zero objects in $\Niq$.

Note also that it follows from Lemma \ref{qabelsi} that

\begin{Proposition}\label{abiso}
Any group isomorphic in $\Niq$ to an abelian group is itself abelian.
\end{Proposition}\qed

\begin{Example}
Let $f:\bZ^3\to\bZ\vee\bZ$ be the map given by
$$
f(l,m,n)=l[x,y]+mx+ny,
$$
where $x$ and $y$ are the generators of $\bZ\vee\bZ$. One then has
$$
\cref f{(l,m,n)}{(l',m',n')}=m'n[x,y],
$$
so that $f$ is a q-map. It is obviously a bijection. However it cannot be an
isomorphism in $\Niq$ because of Proposition \ref{abiso}.

In fact,
$$
\cref{f\1}{l[x,y]+mx+ny}{l'[x,y]+m'x+n'y}=(-m'n,0,0),
$$
so that $f\1$ is quadratic, but not a q-map.
\end{Example}

Let  us point out that there exist nil$_2$-groups isomorphic in $\Niq$ but not
in $\Nil$. We will see such examples below (see Example \ref{died=quat}).

Let us recall that a \emph{weak coproduct} of objects $X_1$ and $X_2$ of some
category is an object $W$ together with morphisms $i_k:X_k\to W$ such that for
any morphisms $f_k:X_k\to Z$ there exists a morphism (not necessarily unique)
$f:W\to Z$ with $f_k=fi_k$, $k=1,2$.

\begin{Lemma}
The category $\Niq$ possesses weak coproducts.
\end{Lemma}

\begin{proof}
We claim that $W=X_1\x X_2$ does the job. Indeed, for any $f_k:X_k\to Z$ put
$f=f_1p_1+f_2p_2$. Then one has
$$
fi_k=(f_1p_1+f_2p_2)i_k=f_1p_1i_k+f_2p_2i_k=f_k.
$$
\end{proof}



\section{The category $\Nil$ as a linear extension}\label{nillinext}

We start with recalling the definition of a linear extension of a small
category by a bifunctor \cite{BW}

\begin{Definition}\label{linearextension}
A \emph{linear extension} of a small category  $\bC$ by a bifunctor $D:\bC\op\x\bC\to\Ab$
$$
0\to D\to\bE\xto P\bC\to0
$$
is a functor $P$  with the following properties: $\bC$  and $\bE$  have the same
objects and $P$  is a full functor which
is the identity on objects.  For each pair of objects $i$ and $j$  the abelian
group $D(i,j)$ acts transitively and effectively on the set $\Hom_\bE(i,j)$. We
write $\alpha + a$ for the action of $a \in D(i,j)$ on $\alpha\in\Hom_\bE(i,j)$.
The action satisfies the linear distributivity law :

$$
(\alpha+a)(\beta+b) = \alpha\beta + P(\alpha)_*b +
P(\beta)^*a.
$$
\end{Definition}

It is known and easy to prove that in any linear extension the functor $q$
reflects isomorphisms and yields a bijection on isomorphism classes of
objects.

%


Our aim is to obtain the category $\Nil$ as a linear extension. To do so we first recall some
classical results on group (co)homology.

\begin{Proposition}\label{H22001}\ 

\begin{itemize}
\item[i)] For a central extension 
\begin{equation}\label{cega}
\bE=\left(0\to A\xto iG\xto pQ\to0\right)
\end{equation}
there is a
well-defined class $\cl\bE\in H^2(Q;A)$ and in this way one obtains a
one-to-one correspondence between the equivalence classes of central extensions
of $Q$ by $A$ and elements of the group $H^2(Q;A)$. If $\bE'$ is also a central
extension of a group $Q'$ by $A'$ and $f:Q\to Q'$, $g:A\to A'$ are group
homomorphisms then $g_*\cl\bE$ and $f^*\cl{\bE'}$ are the same elements in
$H^2(Q;A')$ iff there is a group homomorphism $h:G\to G'$ such that the diagram
$$
\xymatrix{
0\ar[r] &A\ar[r]\ar[d]_g&G\ar[r]\ar[d]_h&Q\ar[d]^f\ar[r]&0\\
0\ar[r] &A'\ar[r]&G'\ar[r]&Q'\ar[r]&0 }
$$
commutes.

\item[ii)] Let $Q$ be a group and $A$ be an abelian group, considered as a $Q$-module via
the trivial action of $Q$ on $A$. Then one has the universal coefficient exact sequence
$$
0\to\Ext(Q\ab,A)\to H^2(Q;A)\xto\mu\Hom(H_2Q,A)\to0.
$$

\item[iii)] For the central extension \eqref{cega} one has the following Ganea
exact sequence
$$
G\ab\ox A\to H_2G\to H_2Q\xto\nu A\to G\ab\to Q\ab\to0,
$$
where $\nu=\mu\cl\bE$, with $\cl\bE$ as in i) above.

\item[iv)] If $B$ is an abelian group then $H_2(B)\cong\Lambda^2B$ and the
homomorphism $\mu\cl\bE:\Lambda^2B\to A$ corresponding to a central extension
$$
\bE=\left(0\to A\xto iG\xto pB\to0\right)
$$
is determined by
$$
i\left(\mu\cl\bE(p(x)\wedge p(y))\right)=[x,y], \ \ x,y\in G.$$
\end{itemize}
\end{Proposition}

\begin{proof}
These results are well known, see for example
\cite{stammbach}.
\end{proof}

The class of the central extension
\begin{equation}\label{grogorc}
0\to[G,G]\to G\to G\ab\to0
\end{equation}
in $H^2(G\ab;[G,G])$ is denoted by ${\mathsf e}(G)$.

\begin{Lemma}\label{1986}
The homomorphism $\mu({\mathsf e}(G)):\Lambda^2(G\ab)\to[G,G]$ is surjective.
\end{Lemma}

\begin{proof}
This follows from iv) of Proposition \ref{H22001} applied to the central
extension \eqref{grogorc}.
\end{proof}

The exact
sequence \eqref{grogorc} is functorial on $G$, meaning that if $f:G\to H$ is a
homomorphism, then one has the following commutative diagram
$$\xymatrix{0\ar[r] &[G,G]\ar[r]^i\ar[d]_{\com f}&G\ar[r]^p\ar[d]_f&
G\ab\ar[d]^{f\ab}\ar[r]&0\\
0\ar[r] &[H,H]\ar[r]^j&H\ar[r]^q&H\ab\ar[r]&0 }
$$
If $g:G\to H$ is another homomorphisms, then we write $f\sim g$ provided
$f\ab=g\ab$ and $\com f=\com g$. It is clear that $f\sim g$ iff there exists
a homomorphism $k:G\ab\to [H,H]$ such that $f-g=jkp$. We can consider the
corresponding quotient category $\Nilsim$. Objects are the same as of $\Nil$. Two
homomorphisms  $f,g:G\to H$ defines the same morphism in $\Nilsim$ provided
$f\sim g$. Comparing with the notion of linear extension of categories
(see Definition \ref{linearextension}) we obtain the following result.

\begin{Theorem}\label{pirveliCrfivi} One has the following linear extension of categories
$$0\to D\to \Nil\to \Nilsim\to 0$$
where the bifunctor
$$D:\left(\Nilsim\right)\op\x\Nilsim\to \Ab$$
is given by
$$D(G,H)=\Hom(G\ab,[H,H]).$$
\end{Theorem}\qed

Our next aim is to describe the quotient category $\Nilsim$ in cohomological
terms. Define the category $\Nilab$ as follows. The objects of $\Nilab$ are
triples $(A,B,e)$, where $A$ and $B$ are abelian groups and $e\in H^2(A,B)$ is
such an elements that $\mu(e):\Lambda^2(A)\to B$ is an epimorphism. A morphism
from $(A,B,e)$ to $(A',B',e')$ is a pair $(f,g)$, where $f:A\to A'$ and $g:B\to
B'$ are homomorphisms such that the equation
$$
f^*(e')=g_*(e)
$$
holds in $H^2(A,B')$. Thus for any $G\in\Nil$ the triple
$$
{\sf ch}(G)=(G\ab,[G,G],{\sf e}(G))
$$
is an object of $\Nilab$. Moreover, if $f:G\to H$ is homomorphism of groups,
then $(f\ab,\com f):{\sf ch}(G)\to {\sf ch}(H)$ is a morphism. In this way
one obtains the functor
$${\sf ch}:\Nilsim\to\Nilab.
$$

\begin{Theorem}\label{19860}
The functor ${\sf ch}:\Nilsim\to\Nilab$ is an equivalence of categories.
\end{Theorem}

\begin{proof} We claim that for any object $(A,B,e)\in\Nilab$
there exist an object $G\in\Nil$ and an isomorphism ${\sf ch}(G)\to(A,B,e)$
in $\Nilab$. Indeed, consider a central extension
$$
0\to B\to G\to A\to 0
$$
corresponding to the element $e$. The exact sequence iii) of Proposition
\ref{H22001} in our case has the following form
$$
H_2G\to H_2A\to B\to G\ab\to A\to 0
$$
Since $H_2A\to B$ is an epimorphism, it follows that $G\ab\cong A$. Therefore
$[G,G]\cong B$ and the claim is proved.

Now, we show that for any  morphism $(f,g):{\sf ch}(G)\to {\sf ch}(G_1)$ in
$\Nilab$, there is a unique homomorphism $h:G\cong G_1$ in $\Nilsim$, such that
${\sf ch}(h)=(f,g)$. Indeed, by definition of morphisms in $\Nilab$ and by
part i) of Proposition \ref{H22001} there exist $h:G\to G_1$ in $\Nil$, such
that the following diagram commutes:
$$
\xymatrix{0\ar[r] &[G,G]\ar[r]^i\ar[d]_g&G\ar[r]^p\ar[d]_h&G\ab\ar[d]^f\ar[r]&0\\
0\ar[r]&[G_1,G_1]\ar[r]^-{i_1}&G_1\ar[r]^{p_1}&G_1\ab\ar[r]&0 }
$$
If $h':G\to G'$ is another such homomorphism, then clearly $h\ab=(h')\ab$ as
well as $\com h=\com{h'}$ and result is proved.
\end{proof}

\section{The category $\Niq$ as a linear extension}\label{niqext}

For $\alpha\in\Hom(G\ab\ox G\ab,[H,H])$ and $f\in\qw(G,H)$ define $f+\alpha\in\qw(G,H)$ to
be the map given by
$$
(f+\alpha)(g)=f(g)+\alpha(\hat g,\hat g).
$$
It is clear that for any $\alpha,\beta\in\Hom(G\ab\ox G\ab,[H,H])$ and $f\in
\qw(G,H)$ one has
$$
f+(\alpha+\beta)=(f+\alpha)+\beta
$$
and therefore the group $\Hom(G\ab\tp G\ab,[H,H])$ acts on the set $\qw(G,H)$.
In particular, this gives the following equivalence relation: for q-maps $f,g\in
\qw(G,H)$ we put $f\sim g$ provided $g=f+\alpha$, for some homomorphism
$\alpha:G\ab\ox G\ab\to[H,H]$.

\begin{Lemma}\label{sab26}\ 

\begin{itemize}
\item[i)] Let $f_1,f_2,g_1,g_2:G\to H$ be q-maps. If $f_1\sim g_1$ and $f_2\sim
g_2$, then
$$
f_1+f_2\sim g_1+g_2.
$$

\item[ii)] Let $f,g:G\to H$ be q-maps. Then
$$
f+g\sim g+f.
$$

\item[iii)] Let $f:G\to H$ and $g_1,g_2:G_1\to G$ be q-maps. Then
$$
f(g_1+g_2)\sim fg_1+fg_2.
$$

\item[iv)] Let $f:G\to H$ and $g:G_1\to G$ be q-maps. Then for any
homomorphisms $\alpha:G\ab\ox G\ab\to[H,H]$ and $\beta:G_1\ab\ox G_1\ab\to
[G,G]$ one has
$$
(f+\alpha)(g+\beta)=fg+f_*\beta +g^*\alpha
$$
where $f_*\beta$ and $g^*\alpha$ are homomorphisms $G_1\ab\tp G_1\ab\to[H,H]$
given by $f_*\beta(x,y)=f(\beta(x,y))$ and
$g^*\alpha(\hat{x},\hat{y})=\alpha(\widehat{g(x)},\widehat{g(y)})$.

\item[v)] Let $f_1,f_2:G\to H$ and $g_1,g_2:G_1\to G$ be q-maps. If  $f_1\sim
f_2$ and $g_1\sim g_2$, then
$$
f_1g_1\sim f_2g_2.
$$

\item[vi)] Let $f,g:G\to H$ be q-maps. If $f\sim g$ then they
induce the same homomorphisms $G\ab\to H\ab$ and $[G,G]\to [H,H]$.
\end{itemize}
\end{Lemma}

\begin{proof} i) We have $g_i=f_i+\alpha_i$,  where $\alpha_i:G\ab\tp G\ab\to
[H,H]$ is a homomorphism $i=1,2$. Since the values of $\alpha_i$ are in the
center, we obtain $g_1+g_2=f_1+f_2+(\alpha_1+\alpha_2)$.

ii) It suffices to observe that $f+g=g+f+\alpha$, where
$\alpha(\hat{x},\hat{y})=[f(x),g(y)]$.

iii) Thanks to equation \eqref{komcr} one has $f(g_1+g_2)=fg_1+fg_2+\alpha$,
where $\alpha(\hat x_1,\hat x_2)=\cref f{g_1(x_1)}{g_2(x_2)}$.

iv) We have
$(f+\alpha)(g+\beta)(x)=f(g(x)+\beta(\hat x,\hat
x))+\alpha(\widehat{g(x)+\beta(\hat x,\hat x)},\widehat{g(x)+\beta(\hat x,\hat
x)})$.
Since the values of $\beta$ lie in the commutator subgroup of $H$ and $\alpha$ is
defined on the abelization, we get
$\alpha(\widehat{g(x)+\beta(\hat x,\hat x)},\widehat{g(x)+\beta(\hat x,\hat
x)})=\alpha(\widehat{g(x)})$. Thus the result follows from iv)
of Lemma \ref{qwadasxvistvisebebi}.

v) This property is an immediate consequence of iv).

vi) By assumption $g=f+\alpha$, for some homomorphism $G\ab\tp G\ab\to [H,H]$.
If $c\in[G,G]$, then $\hat{c}=0$ in $G\ab$, thus $\alpha(\hat{c},\hat{c})=0$
and hence $g(c)=f(c)$. On the other hand, for any $x\in G$ the class of
$\alpha(\hat{x},\hat{x}) $ in $H\ab$ vanishes, hence $f\ab=g\ab$.
\end{proof}

\begin{Corollary} There is a well-defined category $\Niqsim$, with  objects
nil$_2$-groups, and morphisms $\sim$-equivalence classes of q-maps. The
category $\Niqsim$ is an additive category.
\end{Corollary}

\begin{proof} By v) $\Niqsim$ is a well-defined category. By i) and ii) hom's in
$\Niqsim$ are abelian groups. Since $\Niq$ was left distributive, it follows from
iii) that the composition in $\Niqsim$ is distributive. One easily sees that the
product in $\Niq$ remains also a product in $\Niqsim$ and therefore $\Niq$ is an
additive category with products.
\end{proof}

For q-maps $f,g\in\qw(G,H)$ we put $f\approx g$ provided both of them yield the same
homomorphisms $G\ab\to H\ab$ and $[G,G]\to [H,H]$. The corresponding  quotient
category is denoted by $\Niqssim$. By iv) in Lemma \ref{sab26} the quotient functor
$\Niq\to\Niqssim$ factors trough $\Niqsim$.

For nil$_2$-groups $G$ and $H$ we let $D^\sim(G,H)$ be the quotient of
$\Hom(G\ab\tp G\ab,[H,H])$ by the subgroup spanned by such $\alpha\in
\Hom(G\ab\tp G\ab,[H,H])$ that $\alpha(\hat{x},\hat{x})=0$ for all $x\in G$. In
this way one obtains a bifunctor
$D^\sim:\left(\Niqssim\right)\op\x\Niqssim\to\Ab$.  We also
need another bifunctor $D^\approx:\left(\Niqssim\right)\op\x\Niqssim\to\Ab$ given by
$D^\approx(G,H)=\qu(G\ab,[H,H])$. There is a natural transformation
$\rho:D^\sim\to D^\approx$, which takes $\alpha:G\ab\tp G\ab\to [H,H]$ to the quadratic map
$\rho(\alpha):G\ab\to [H,H]$ given by
$\rho(\alpha)(\hat{x})=\alpha(\hat{x},\hat{x})$. It follows from the definition
of $D^\sim$, that $\rho$ is a monomorphism. We define $\tilde{D}:=\cok(\rho)$. Using
the quotient functors $\Niq\onto\Niqsim\onto\Niqssim$ one  considers
$D^\sim$, $D^\approx$ also as bifunctors on $\Niqsim$, or $\Niq$.

\begin{Proposition} One has the following commutative diagram of linear extensions:
$$
\xymat{
&
&0\ar[d]
&0\ar[d]\\
0\ar[r]
&D^\sim\ar@{=}[d]\ar[r]^\rho
&D^\approx\ar[r]\ar[d]
&\tilde D\ar[r]\ar[d]&0\\
0\ar[r]
&D^\sim\ar[r]
&\Niq\ar[r]\ar[d]
&\Niqsim\ar[r]\ar[d]&0\\
&
&\Niqssim\ar@{=}[r]\ar[d]
&\Niqssim\ar[d]\\
&
&0
&0
&
}
$$
In particular $\Niqssim$ is also an additive category and the quotient functors
$\Niq\to\Niqsim\to\Niqssim$ reflect isomorphisms and yield bijections on
isomorphism classes of objects.
\end{Proposition}

\begin{proof}  The operation $\qw(G,H)\x\Hom(G\ab\tp G\ab,[H,H])\to\qw(G,H)$ given by
$(f,\alpha)\mapsto f+\alpha$ yields the action of $D^\sim$ on the category $\Niq$
and by the property iv) one obtains a linear extension of categories
$$
0\to D^\sim\to\Niq\to\Niqsim\to 0.
$$
By Proposition \ref{2x3} for q-maps $f,g:G\to H$ one has $f\approx g$ iff there
is a quadratic map $h:G\ab\to [H,H]$ such that $f-g$ factors trough $h$. This
shows that
$$
0\to D^\approx\to\Niq\to\Niqssim\to 0
$$
is a linear extension of categories. The rest follows from the properties of
linear extensions.
\end{proof}

\begin{Remark} For an abelian group $A$ the group $A\tp A$ has a canonical
involution $(a\tp b)^\sigma=b\tp a$. We put $\tilde{\Gamma}^2(A):=\{x\in A\tp
A\mid x^\sigma=x\}$. Then one has an exact sequence
$$0\to \tilde{\Gamma}^2(A)\to A\tp A\to \Lambda^2(A)\to 0$$
One easily sees that the class of abelian groups for which this sequence splits
is closed under direct sums and contains all cyclic groups (and all uniquely
2-divisible groups). In particular the sequence
splits, provided $A$ is a direct sum of cyclic groups. The exact sequence
yields the following exact sequence
$$0\to \Hom(\Lambda^2(A),B)\to \Hom(A\tp A,B)\xto{\xi_{A,B}} \Hom(\tilde{\Gamma}^2(A),B)$$
for all abelian group $B$. It follows from the definition that $D^\sim(G,H)\cong
\im(\xi_{G\ab,[H,H]})$. In particular, if $G\ab$ is  a direct sum of cyclic
groups, then $D^\sim(G,H)=\Hom (\tilde{\Gamma}^2(G\ab), [H,H]).$
\end{Remark}

\begin{Definition}\label{h2b}
For abelian groups $A$ and $B$ we denote by $H_b^2(A,B)$ be the subgroup of
$H^2(A,B)$ generated by bilinear 2-cocycles. Thus by definition one has the
following exact sequence
$$
\qu(A,B)\xto{\cref?--}\Hom(A\tp A,B)\to H_b^2(A,B)\to 0
$$
where the first map assigns to a quadratic map $f$ its cross-effect
$\cref f--$.

We now define the category $\Niqab$, which has the same objects as the
category $\Nilab$. Thus objects are triples $(A,B,e)$ where $A$ and $B$ are
abelian groups, and $e\in H^2(A,B)$ is such an element that $\mu(e):\Lambda^2(A)\to
B$ is an epimorphism. A morphism from $(A,B,e)$ to $(A',B',e')$ in $\Niqab$
is a pair $(f,g)$, where $f:A\to A'$ and $g:B\to B'$ are homomorphisms such
that
$$
f^*(e')-g_*(e)\in H_b^2(A,B').
$$
\end{Definition}

\begin{Theorem}\label{niq''}
The functor ${\sf ch}:\Nilsim\to\Nilab$ has a canonical extension
$$
{\sf ch}:\Niqssim\to\Niqab
$$
which is an equivalence of categories.
\end{Theorem}

\begin{proof} On objects one puts
$$
{\sf ch}(G)=(G\ab,[G,G],{\sf e}(G)).
$$
If $f:G\to H$ is a q-map, then one puts
$$
{\sf ch}(G\xto{f} H)=( f\ab,\com f).
$$
We claim that one has
$$
\left(f\ab\right)^*({\sf e}(H))-{\com f}_*({\sf e}(G))\in H_b^2(G\ab,[H,H]).
$$
Let $\alpha$ (resp. $\beta$) be a 2-cocycle representing the class ${\sf e}(G)$
(resp. ${\sf e}(H)$). Thus $G=G\ab\x [G,G]$ (resp. $H=H\ab\x[H,H]$) as a set,
with the following group structure $(a,u)+(b,v)=(a+b,\alpha(a,b)+u+v)$, where
$a,b\in G\ab$ and $u,v\in [G,G]$ (resp. $(c,x)+(d,y)=(c+d,\beta(c,d)+x+y)$,
$c,d\in H\ab$, $x,y\in [H,H]$). Any q-map $f:G\to H$ has the form
$f(a,u)=(f\ab(a),\com f(u)+\gamma(a))$, where $f\ab:G\ab\to H\ab$ and
$\com f:[G,G]\to [H,H]$ are homomorphisms, while $\gamma:G\ab\to[H,H]$ is a
map. One has
\begin{align*}
f((a,u)+(b,v))
=&f(a+b,\alpha(a,b)+u+v)\\
=&(f\ab(a)+f\ab(b),\gamma(a+b)+\com f(\alpha(a,b))+\com f(u)+\com f(v)).
\end{align*}
On the other hand we have $f((a,u)+(b,v))=f((a,u))+f((b,v))+((a,u),(b,v))_f$.
Since the cross-effect of $f$ factors through $\delta:G\ab\ox G\ab\to [H,H]$ we
obtain
\begin{align*}
&f((a,u)+(b,v))\\
&=(f\ab(a),\com f(u)+\gamma(a))+(f\ab(b),\com f(v)+\gamma(b))+(0,\delta(a,b))\\
&=(f\ab(a)+f\ab(b),\beta(f\ab(a),f\ab(b))+\gamma(a)+\gamma(b)+\delta(a,b)+\com
f(u)+\com f(v)).
\end{align*}
Comparing these expressions we obtain
$$
\gamma(a+b)+\com f(\alpha(a,b))=\beta(f\ab(a),f\ab(b))+\gamma(a)+\gamma(b)+\delta(a,b).
$$
Thus the class $\left(f\ab\right)^*({\sf e}(H))-{\com f}_*({\sf e}(G))$ in the group
$H_b^2(G\ab,[H,H])$ coincides with the class of $-\delta$ and the claim is
proved. It follows that $\sf ch$ is a well-defined functor $\Niq\to
\Niqab$, which obviously factors trough the category $\Niqssim$. By our
construction and by definition of $\Niqssim$ the induced map
$$
\Hom_\Niqssim(G,H)\to \Hom_\Niqab({\sf ch}(G),{\sf ch}(H))
$$
is an injection. Let us show that this map is surjective as well. Take any
morphism $(g,h):{\sf ch}(G)\to {\sf ch}(H)$ in $\Niqab$. Then $g:G\ab\to
H\ab$ and $h:[H,H]\to [G,G]$ are homomorphisms such that
$$
h(\alpha(a,b))-\beta(g(a),g(b))=-\gamma(a+b)+\gamma(a)+\gamma(b)+\delta(a,b)
$$
where $\delta:G\ab\tp G\ab\to [H,H]$ is a homomorphism, $\gamma:G\ab\to [H,H]$
is a map, while $\alpha$ and $\beta$ are as above. Define the map $f:G\to H$ 
via
$f(a,u)=(g(a),\gamma(a)+h(u))$. Then one has
$$((a,u),(b,v))_f=\delta(a,b).$$
Thus $f$ is a q-map with $f\ab=g$ and $\com f=h$. Therefore $\sf ch$ is full
and faithful. By Theorem \ref{19860} the functor $\sf ch$ is surjective on
isomorphism classes of objects and the result follows.
\end{proof}

\section{$\mathrm q$-split groups}\label{qsplit}
We start with the following definitions.

\begin{Definition}\label{similar}
Call nil$_2$-groups \emph{similar} if they have isomorphic abelianizations and
isomorphic commutator subgroups.
\end{Definition}

\begin{Definition}
Call a nil$_2$-group $G$ \emph{q-split} if the quotient map $G\onto G\ab$ has a
quadratic section. It is easy to see that this section is then a q-map.
\end{Definition}
\begin{Lemma}\label{qsprod}
The class of q-split groups contains all abelian groups and is closed under
products and coproducts.
\end{Lemma}

\begin{proof}
For products and abelian groups this is obvious. For coproducts, note that
 the central extension
$$
0\to G_1\ab\ox G_2\ab\to G_1\vee G_2\to G_1\x G_2\to0
$$
has a quadratic section $s$ given by $s(g_1,g_2)=(0,g_1,g_2)$. One easily
checks that
$$
\cref s{(x_1,x_2)}{(y_1,y_2)}=(y_1\ox x_2,0,0)=[(0,y_1,0),(0,0,x_2)].
$$
Thus $s$ is a q-map. Since $(G_1\vee G_2)\ab=(G_1\x G_2)\ab=G_1\ab\x G_2\ab$
we see that for any quadratic sections $s_i:G_i\ab\to G_i$ of the natural
projections $G_i\onto G_i\ab$, $i=1,2$, the composite $s\circ(s_1\x
s_2):(G_1\vee G_2)\ab\to G_1\vee G_2$ is a section. Since $s$, $s_1$, $s_2$
are q-maps, $s\circ(s_1\x s_2)$ is also a q-map and  the result follows.
\end{proof}

\begin{Example}\label{died=quat}
It follows that the dihedral group $D_4\cong \bZ/2\bZ\vee \bZ/2\bZ$ of order
$8$ is q-split. Let us show that the quaternion group
$Q_8=\brk{\tau,\omega\mid 2\tau=2\omega,\omega+\tau-\omega=-\tau}$ of order
$8$ is also q-split. Observe that  $\tau$ and $\omega$ are of order $4$ and
$[\omega,\tau]=2\tau$. So one has $Q_8\ab\cong \bZ/2\bZ\x \bZ/2\bZ\cong
D_4\ab$ and $[Q_8,Q_8]\cong \bZ/2\bZ\cong [D_4,D_4]$. One easily checks that
the map $s:Q_8\ab\to Q_8$ given by $s(0)= 0$, $s(\hat{\omega})=\omega$,
$s(\hat{\tau})= \tau$, $s(\hat{\omega}+\hat{\tau})= \omega+\tau$ is a
quadratic section of $Q_8\to Q_8\ab$. It follows from Corollary \ref{qspliso}
below  that $Q_8$ and $D_4$ are isomorphic in $\Niq$.
\end{Example}

\begin{Lemma}\label{einb}
A nil$_2$-group $G$ is q-split iff the class ${\sf e}(G)\in H^2(G\ab,[G,G])$
belongs to the subgroup $H^2_b(G\ab,[G,G])$.
\end{Lemma}

\begin{proof}
Let $u:G\ab\to G$ be a quadratic section. Then the class ${\sf e}(G)$ can be
represented by the  cocycle $(a_1,a_2)\mapsto u(a_1)+u(a_2)-u(a_1+a_2)$ which
is bilinear and therefore lies in $H^2_b(G\ab,[G,G])$. Conversely, if the
class ${\sf e}(G)\in H^2(G\ab,[G,G])$ is represented by a bilinear map
$f:G\ab\x G\ab\to [G,G]$, then $G$ is isomorphic to the set $G\ab\x [G,G]$
with group structure defined by
$(a_1,b_1)+(a_2,b_2)=(a_1+a_2,b_1+b_2+f(a_1,a_2))$, and the projection of the
latter to $G\ab$ has a quadratic section given by $a\mapsto(a,0)$.
\end{proof}

We denote by $\SpN$ and $\SpNN$ the full subcategories of, respectively,
$\Niq$ and $\Niqssim$ with objects all q-split groups. They are related via
the following linear extension:
$$
0\to D^\approx\to \SpN\to\SpNN\to 0
$$
and in particular they have the same isoclasses of objects. According to Theorem \ref{niq''} and
Lemma \ref{einb} the category  $\SpNN$ is equivalent to the category $\Splab$,
which is the full subcategory of the category $\Nilab$ on those objects
$(A,B,e)$ of $\Nilab$ satisfying $e\in H^2_b(A,B)$.  Let us observe that 
$$\Hom_{\Nilab}((A,B,e), (A',B',e'))=\Hom(A,A')\x \Hom(B,B')$$
because the compatibility condition required in the definition of morphisms in $\Nilab$
holds authomatically in $\Splab$.

We now consider another category $\Spl$, which is a full subcategory of the
product category $\Ab\x\Ab$. Objects of the category $\Spl$ are pairs of
abelian groups $(A,B)$ for which there exists a homomorphism $f:A\tp A\to B$
such that $f^a:\Lambda ^2(A)\to B$ is an epimorphism, where
$f^a(x,y):=f(x,y)-f(y,x)$.

\begin{Theorem}\label{qspliso}
The categories $\Spl$ and $\Nilab$ are equivalent. Thus, two q-split groups 
 are isomorphic in $\Niq$ iff they are similar.
\end{Theorem}

\begin{proof} Take any object $(A,B)$ of $\Spl$ and choose $f:A\tp A\to B$ for which $f^a$ is
an epimorphism. Let $e_f\in H^2_b(A,B)$ be the class correspobding to $f$. Then 
$(A,B,e_f)\in \Nilab$. Then $(A,B,f)\mapsto (A,B,e_f)$ yields expected equivalence of categories.
\end{proof}

\begin{Remark}
One easily sees that the class $\bsf S$ of abelian groups for which the
natural short exact sequence $0\to\Lambda^2(A)\to A\tp A\to\sym^2(A)\to 0$
splits contains all cyclic groups, all uniquely 2-divisible groups and is
closed under direct sums. In particular any finitely generated abelian group
lies in $\bsf S$. If $A\in\bsf S$ then for any homomorphism $g:\Lambda^2(A)\to
B$ there exists a homomorphism $f:A\ox A\to B$ such that $f^a=g$. It follows
that a pair of abelian groups $(A,B)$ with $A\in\bsf S$ belongs to $\Spl$ iff
there exists an epimorphism $\Lambda^2(A)\to B$.
\end{Remark}

\begin{Proposition}
For abelian groups $A$, $B$ there is a commutative diagram with exact rows
\begin{equation}\label{qdia}
\xymat{
0\ar[r]
&\Hom(\sym^2A,B)\ar[r]\ar@{-->}[d]
&\Hom(A\ox A,B)\ar[r]\ar[d]^\alpha
&\Hom(\Lambda^2A,B)\ar@{=}[d]\\
0\ar[r]
&\Ext(A,B)\ar[r]
&H^2(A;B)\ar[r]^-c
&\Hom(\Lambda^2A,B)\ar[r]
&0}
\end{equation}
where the image of the homomorphism $\alpha$ is equal to the subgroup
$H^2_b(A,B)$ from \ref{h2b}.
\end{Proposition}

\begin{proof}
For any abelian group $A$ one has a short exact sequence
$$
0\to\Lambda^2A\to A\ox A\to\sym^2A\to0,
$$
We place in the upper row of \eqref{qdia} the sequence induced by this short
exact sequence. The lower row is the universal coefficient exact sequence, and
the map $\alpha$ is given by considering a bilinear map as a 2-cocycle.
The rest is obvious.
\end{proof}

\begin{Remark}
The arrow on the upper right of \eqref{qdia} is surjective if $A$ is either
uniquely 2-divisible or is a direct sum of cyclic groups. Indeed in these
cases the aforementioned short exact sequence splits.
\end{Remark}

\begin{Proposition}
If 2 is invertible in $B$ then the above homomorphism
$\Hom(\sym^2A,B)\to\Ext(A,B)$ is zero.
\end{Proposition}

\begin{proof}
For a homomorphism $f:\sym^2A\to B$ the class $\alpha(f)$ is represented by the
cocycle $(x,y)\mapsto f(xy)$. This cocycle is the coboundary of the cochain
$g:A\to B$ given by $a\mapsto\2f(a^2)$.
\end{proof}

\begin{Lemma}\label{nonqsp}
Let $A$ be any abelian group and let $B$ be a uniquely 2-divisible group. Then
for any $x\in H^2_b(A;B)$ and any $0\ne a\in\Ext(A,B)$ one has $x+a\notin H^2_b(A,B)$. 
\end{Lemma}

\begin{proof} 
Otherwise one would have $a\in\im\alpha$, which contradicts the previous lemma.
\end{proof}

\begin{Remark}
It follows that for any object $(A,B,x)$ of $\Niqab$ and any $a\in\Ext(A,B)$,
also $(A,B,x+a)$ gives an object of $\Niqab$, since in the universal
coefficient exact sequence in \eqref{qdia} one has $c(x+a)=c(x)$. In
particular, if $(A,B,x)$ with uniquely 2-divisible $B$ lies in the
subcategory $\Splab$ and $a\ne0$, then $(A,B,x+a)\in\Niqab$ cannot belong to
$\Splab$, since by Lemma \ref{nonqsp}, $x+a\notin H^2_b(A;B)$. Thus not all
nil$_2$-groups are q-split. Some explicit examples of non-q-split groups
follow.
\end{Remark}

\begin{Example}
For each prime $p$ consider the semidirect product
$\bZ/p^2\bZ\ltimes\bZ/p\bZ$, where the generator of $\bZ/p\bZ$ acts on
$\bZ/p^2\bZ$ via multiplication by $p+1$. This group is similar in the sense
of Definition \ref{similar} to $\bZ/p\bZ\vee\bZ/p\bZ$ (both have
abelianizations isomorphic to $(\bZ/p\bZ)^2$ and commutator subgroups
isomorphic to $\bZ/p\bZ$). For $p=2$ these groups are in fact isomorphic;
however for odd $p$ they are not, since the former has exponent $p^2$ and the
latter has exponent $p$. Thus in the diagram \eqref{qdia} for
$A=(\bZ/p\bZ)^2$ and $B=\bZ/p\bZ$, classes of $\bZ/p\bZ\vee\bZ/p\bZ$ and
$\bZ/p^2\bZ\ltimes\bZ/p\bZ$ in $H^2((\bZ/p\bZ)^2;\bZ/p\bZ)$ are not equal. On
the other hand one can choose isomorphisms of their commutator subgroups with
$\bZ/p\bZ$ in a way which makes obvious that these classes have the same
image under the homomorphism $c$ defined in \eqref{qdia} above, hence they
differ by a nonzero element of $\Ext(\bZ/p\bZ,\bZ/p\bZ)$. But
$\bZ/p\bZ\vee\bZ/p\bZ$ is q-split by Lemma \ref{qsprod}, hence its class is
in the image of the homomorphism $\alpha$ from \eqref{qdia}. Then by Lemma
\ref{nonqsp} we conclude that $\bZ/p^2\bZ\ltimes\bZ/p\bZ$ is not q-split. In
particular, the above similar groups are also not isomorphic in $\Niq$.
\end{Example}


\section{$\mathrm q$-maps for uniquely 2-divisible groups}\label{liecase}

Let us recall the relevant part of the classical Maltsev correspondence
between nilpotent groups and Lie algebras. In the nil$_2$ case it amounts to
an isomorphism of categories from the category of nil$_2$ Lie algebras over
$\bZ[\2]$, i.~e. Lie algebras with $[L,[L,L]]=0$ to the category of uniquely
2-divisible nil$_2$-groups. In what follows, all Lie algebras are understood
to be of the above kind, i.~e. class two nilpotent Lie $\bZ[\2]$-algebras. Let us
denote by $\Nil(\bZ[\2])$ the category of these algebras and their
homomorphisms. Moreover we will denote by $\Nil^\2$ the category of uniquely
2-divisible nil$_2$-groups.

One defines an isomorphism of categories
$$
\exp:\Nil(\bZ[\2])\to\Nil^\2
$$
by declaring, for an algebra $L\in\Nil(\bZ[\2])$, $\exp(L)$ to be the set $L$
equipped with the operation
$$
a\oplus b=a+b+\2[a,b].
$$
This is a group, with zero element 0 and inverse of an element $a$ given by
$-a$. Moreover the commutator with respect to this group structure coincides
with the Lie bracket, so that for any $L$ one has $[\exp(L),\exp(L)]=[L,L]$
and $\exp(L)\ab=L\ab$, where $L\ab=L/[L,L]$ is the abelianization of the Lie algebra
$L$.

Now clearly any Lie algebra homomorphism is also a homomorphism with respect
to $\oplus$. Moreover we have $a\oplus c=a+c$ for $c\in[L,L]$, so that for
any $a,b\in L$
$$
a+b=a\oplus b\oplus\2[b,a].
$$
It follows that also conversely, a map which is a homomorphism with respect to
$\oplus$ is a Lie algebra homomorphism, so that $\exp$ is an isomorphism of
categories, with the inverse isomorphism $\log$ defined as follows: for a
uniquely 2-divisible nil$_2$-group $G$ the Lie algebra $\log(G)$ is the set
$G$ equipped with the addition as above and with the bracket equal to the
commutator map.

Our aim in this section is to prove

\begin{Theorem}\label{bolo}
Two uniquely 2-divisible nil$_2$ groups $G$, $G'$ are isomorphic as objects of $\Niq$ if
and only if there exists an isomorphism of abelian groups
$g:\log(G)\to\log(G')$ such that $g[G,G]=[G',G']$.
\end{Theorem}

For the proof we must define an analog of the category $\Niq$ from Section \ref{niq}
for Lie algebras. For this, we first define

\begin{Definition}\label{qlie}
A map $f:L\to L'$ between Lie algebras is called a
q-map if it is a quadratic map between the underlying abelian groups
and moreover for any $a,b\in L$ and any $c\in[L,L]$ one has
$\cref fab\in[L',L']$, $f(a+c)=f(a)+f(c)$ and $f(c)\in[L',L']$.

Moreover we consider the following category
$\Niq(\bZ[\2])\supset\Nil(\bZ[\2])$ with the same objects as $\Nil(\bZ[\2])$.
A morphism $L\to L'$ in $\Niq(\bZ[\2])$ is a q-map in the sense just defined.
\end{Definition}

The key observation is then

\begin{Theorem}\label{niso}
The functor $\exp$ extends to an isomorphism of categories
$$
\Niq^\2\simeq\Niq(\bZ[\2]).
$$
\end{Theorem}

This theorem follows immediately from the following

\begin{Proposition}\label{niqlie}
Let $f:L\to L'$ be a map between Lie algebras. Then the following assertions are equivalent:
\begin{itemize}
\item[i)]\label{isqlie}
$f$ is a q-map in the sense of \ref{qlie};
\item[ii)]\label{isqgr}
 $f$ is a q-map when considered as a map $\exp(L)\to\exp(L')$;
\item[iii)]\label{gh}
there exists a linear map $g:L\to L'$ with $g[L,L]\subseteq[L',L']$ and a
symmetric bilinear map $h:L\ab\x L\ab\to[L',L']$ such that one has
$$
f(a)=g(a)+\2h(\hat a,\hat a)
$$
for all $a\in L$.
\end{itemize}
\end{Proposition}

\begin{proof}

ii) $\iff$ iii):

Let $\cref fab^+$, $\cref fab^\oplus$ denote the cross-effect of $f$ with
respect to the corresponding operations. Thus $f$ is a q-map when considered
as a map $\exp(L)\to\exp(L')$ iff $\cref fab^\oplus$ is bilinear and lands in
$[L',L']$. In that case we have
$$
f(a+b)=f(a\oplus b\oplus\2[b,a])=f(a\oplus b)\oplus f(\2[b,a])=fa\oplus
fb\oplus\cref fab^\oplus\oplus\2f[b,a]
$$
and
$$
fa+fb=fa\oplus fb\oplus\2[fb,fa],
$$
hence
\begin{align*}
&\cref fab^+=-(fa+fb)+f(a+b)=-(fa+fb)\oplus f(a+b)\oplus\2[f(a+b),-(fa+fb)]\\
=&\2[fa,fb]\oplus\cref fab^\oplus\oplus\2f[b,a]\oplus\2[fa\oplus
fb\oplus\cref fab^\oplus\oplus\2f[b,a],-(fa\oplus fb\oplus\2[fb,fa])]\\
=&\2[fa,fb]\oplus\cref fab^\oplus\oplus\2f[b,a].
\end{align*}
The latter expression is then symmetric since it is the cross-effect of some
map with respect to the commutative operation $+$. It is bilinear with respect
to $\oplus$ and satisfies
$$
\cref f{a\oplus c}b^+=\cref fa{b\oplus c}^+=\cref fab^+
$$
for any $c\in[L,L]$ and any $a,b\in L$. Hence it is also bilinear with respect
to $+$ and defining
$$
h(\acl a,\acl b)=\cref fab^+
$$
gives a well-defined symmetric bilinear map $h:L\ab\x L\ab\to[L',L']$. Then
the map $g:L\to L'$ given by
$$
g(a)=f(a)-\2h(\hat a,\hat a)=f(a)-\2\cref faa^\oplus
$$
carries $[L,L]$ to $[L',L']$. Moreover this map is linear since
$\cref faa^\oplus=\cref faa^+$ for any $a\in L$, so that
\begin{align*}
g(a+b)=&f(a+b)-\2\cref f{a+b}{a+b}^+\\
=&fa+fb+\cref fab-\2\cref faa^+-\2\cref fbb^+-\2\cref fab^+-\2\cref fba^+\\
=&g(a)+g(b).
\end{align*}

Conversely, given $g$ and $h$ as in iii), we compute
\begin{align*}
&\cref fab^\oplus=-(fa\oplus fb)\oplus f(a\oplus b)\\
=&-(fa+fb+\2[fa,fb])+f(a+b+\2[a,b])+\2[fa+fb+\2[fa,fb],f(a+b+\2[a,b])]\\
=&-(ga+\2h(\hat a,\hat a)+gb+\2h(\hat b,\hat b)+\2[ga,gb])\\
&+ga+gb+\2g[a,b]+\2h(\acl a+\acl b,\acl a+\acl b)+\2[-(ga+gb),ga+gb]\\
=&-\2[ga,gb]+\2g[a,b]+h(\acl a,\acl b)
\end{align*}
which lies in $[L',L']$ and is bilinear, so indeed $f$ is a q-map.

i) $\iff$ iii):

Obviously any $f$ satisfying iii) is quadratic. Moreover,
a map $f$ between $\bZ[\2]$-modules is quadratic if and only if it has the
form
$$
f(a)=g(a)+\2h(a,a)
$$
for unique linear map $g$ and bilinear symmetric map $h$. One just takes
$g(a)=2f(a)-\2f(2a)$ and $h(a,b)=f(a+b)-f(a)-f(b)$. Then it is easy to check
that a quadratic map is a q-map of Lie algebras if and only if the
corresponding $g$ and $h$ satisfy conditions in iii).
\end{proof}

This enables us to obtain an extension to the q-world of the above classical
Maltsev equivalence, by identifying the full subcategory $\Niq^\2\subset\Niq$
on the uniquely 2-divisible nil$_2$-groups with the following category
defined in terms of Lie $\bZ[\2]$-algebras.

Moreover in this situation \ref{niq''} admits a strengthening. To formulate it we
will need some more categories.

\begin{Definition}\label{bevri}
Let $\Niq_0(\bZ[\2])\subset\Niq(\bZ[\2])$ be the subcategory with the same objects and
those morphisms which are actually linear. That is, a morphism from $L$ to
$L'$ in $\Niq_0(\bZ[\2])$ is an abelian group homomorphism $g:L\to L'$ with
$g[L,L]\subseteq[L',L']$.

Moreover let $\widetilde\Niq_0(k)$ be the quotient category of $\Niq_0(k)$
obtained by identifying those $g_1,g_2:L\to L'$ for which
$\restr{g_1}{[L,L]}=\restr{g_2}{[L,L]}$ and $g_1\ab=g_2\ab:L\ab\to{L'}\ab$.
\end{Definition}
 
We then have

\begin{Proposition}\label{nlext}
There are linear extensions
$$
0\to\Ph\to\Niq(\bZ[\2])\xto q\Niq_0(\bZ[\2])\to0
$$
and
$$
0\to\tilde\Ph\to\Niq_0(\bZ[\2])\xto{\tilde q}\widetilde\Niq_0(\bZ[\2])\to0
$$
defined as follows. The functor $q:\Niq(\bZ[\2])\to\Niq_0(\bZ[\2])$ is identity on
objects and given on morphisms via
$$
q(f)(a)=2f(a)-\2f(2a).
$$
The bifunctor $\Ph:\Niq_0(\bZ[\2])\op\x\Niq_0(\bZ[\2])\to\Ab$ is given by
$$
\Ph(L,L')=\Hom(\sym^2(L\ab),[L',L']).
$$
The  functor $\tilde q$ is the canonical quotient functor, and $\tilde\Ph$ is
given by
$$
\tilde\Ph(L,L')=\Hom(L\ab,[L',L']).
$$
Moreover the categories $\Niq_0(k)$ and $\widetilde\Niq_0(k)$ are both
additive, and the functor $q$ has a section given by the embedding.
\end{Proposition}

\begin{proof}
Additivity of $\Niq_0(\bZ[\2])$ follows from the obvious fact that for any morphisms
$g_1,g_2:L\to L'$ in $\Niq_0(\bZ[\2])$ the maps $g_1\pm g_2$ are morphisms of
$\Niq_0(\bZ[\2])$ too.

The rest is clear in view of the above considerations. Indeed
we can replace a morphism $f:L\to L'$ in $\Niq(\bZ[\2])$ by a
pair $(g,h)$ as in iii) of Proposition \ref{niqlie}. Under this identification
the functor $q$ becomes the projection sending $(g,h)$ to $g$ and the first
linear extension becomes obvious. The second one is straightforward.
\end{proof}

\begin{Definition}
Let $\Niqab(\bZ[\2])$ denote the following category.
Objects of $\Niqab(\bZ[\2])$ are short exact sequences
$$
0\to B\to E\to A\to0
$$
of $\bZ[\2]$-modules such that there exists a surjective homomorphism
$\pi:\Lambda^2(A)\onto B$. A morphism from $0\to B\to E\to A\to0$ to $0\to B'\to
E'\to A'\to0$ is a pair $(\alpha:A\to A',\beta:B\to B')$ of
homomorphisms such that there exists $\eps:E\to E'$ making the diagram
$$
\xymat{
0\ar[r]&B\ar[r]\ar[d]_\beta&E\ar[r]\ar@{-->}[d]_\eps&A\ar[r]\ar[d]^\alpha\ar[r]&0\\
0\ar[r]&B'\ar[r]&E'\ar[r]&A\ar[r]&0
}
$$
commute. We do not make $\pi$ or $\eps$ part of the structure, in particular
$\pi$ is not required to be compatible with $\alpha$ and $\beta$ in any way.
\end{Definition}

Note that $\Niqab(\bZ[\2])$ is an additive category, since for any $A$, $A'$ there
are surjective homomorphisms $\Lambda^2(A\oplus
A')\onto\Lambda^2(A)\oplus\Lambda^2(A')$ and moreover for any morphism
$(\alpha,\beta)$ in $\Niqab(\bZ[\2])$ the pair $(-\alpha,-\beta)$ is also a
morphism.

There is a functor $r:\widetilde\Niq_0(\bZ[\2])\to\Niqab(\bZ[\2])$ sending a Lie algebra $L$ to the
short exact sequence
$$
0\to[L,L]\to L\to L\ab\to0
$$
and the morphism $[g]:L\to L'$ to the pair $(g\ab,\restr g{[L,L]})$, where
$[g]$ denotes the equivalence class of $g$ and $g\ab:L\ab\to{L'}\ab$ is the
homomorphism induced by $g$ which exists since $g[L,L]\subseteq[L',L']$.

\begin{Proposition}\label{nikab}
The above functor $r$ yields an equivalence of categories
$$
\widetilde\Niq_0(\bZ[\2])\simeq\Niqab(\bZ[\2]).
$$
\end{Proposition}

\begin{proof}
First, $r$ is surjective on objects, since for any object $0\to B\to
E\to A\to0$ of $\Niqab(k)$ any surjective homomorphism $\Lambda^2(A)\onto B$
determines a bracket
$$
[,]:\Lambda^2(E)\onto\Lambda^2(A)\onto B\into E
$$
on $E$ which turns it into a nil$_2$ Lie algebra with $[E,E]=B$ and $E\ab=A$.

Next, $r$ is full since by definition a morphism from the object
$0\to[L,L]\to L\to L\ab\to0$ to the object $0\to[L',L']\to L'\to{L'}\ab\to0$
in $\Niqab(\bZ[\2])$  is by definition a pair of linear maps
$\beta:[L,L]\to[L',L']$, $\alpha:L\ab\to{L'}\ab$ for which there exists a
linear map $g:L\to L'$ fitting in the appropriate diagram, which means that
$\beta=\restr g{[L,L]}$ and $\alpha=g\ab$.

Finally $r$ is faithful since for $g_1,g_2:L\to L'$ one has $r[g_1]=r[g_2]$ if
and only if $g_1$ and $g_2$ are equivalent in the sense of \ref{bevri}, i.~e.
if and only if $[g_1]=[g_2]$.
\end{proof}

We can now finish the proof of our theorem.

\begin{proof}[Proof of \ref{bolo}]
There is a chain of functors
$$
\Niq^\2\xto{\eqref{niso}}
\Niq(\bZ[\2])\xto{\eqref{nlext}}
\Niq_0(\bZ[\2])\xto{\eqref{nlext}}
\widetilde\Niq_0(\bZ[\2])\xto{\eqref{nikab}}
\Niqab(\bZ[\2])
$$
each of which is either an equivalence or a linear extension. The statement of
\ref{bolo} is that objects on the left are isomorphic if and only if their
images under the composite functor are. This is clear since any linear
extension reflects isomorphy of objects.
\end{proof}

\section{A cohomological obstruction to $\mathrm q$-splitting}\label{obstr}
We start
with recalling the definition of
the nonabelian cohomology.
Let $G^*$ be a cosimplicial group. One denotes by $\pi^0(G^*)$ the subgroup
of $G^0$ consisting of elements $x\in G^0$ such that $d^0(x)=d^1(x)$.
Moreover, one defines the pointed set $\pi^1(G^*)$ as the quotient of the
pointed set
$$
Z^1(G^*)=\set{y\in G^1\mid d^1(y)=d^0(y)+d^2(y)}
$$
by the following equivalence relation: $y\sim z$, $y,z\in Z^1(G^*)$  iff
there exists $x\in G^0$ such that $z=-d^0x+y+d^1x$. If $G^*$ is abelian then
one defines $\pi^*(G^*)$ in all dimensions using the homology of the associated
cochain complex $(G^*,d=\sum(-1)^id^i)$.  In particular $\pi^i(G^*)$ is an
abelian group, $i\geq 0$, provided $G^*$ is abelian cosimplicial group. The
following result is well known.

\begin{Lemma}\label{exactnonabelianhomotopy} Let
$$0\to A^*\to G^*\to B^*\to 0$$
be a short exact sequence of cosimplicial groups. Then one has the exact
sequence of pointed sets:
$$0\to \pi^0(A^*)\to \pi^0(G^*)\to \pi^0(B^*)\to \pi^1(A^*)\to
\pi^1(G^*)\to \pi^1(B^*)$$ Moreover, if $A^*$ is abelian, then the connecting
map $\pi^0(B^*)\to \pi^1(A^*)$ is a homomorphism.
\end{Lemma}


We also need the following

\begin{Lemma}\label{quqwzustia} Let $G$ and $H$ be nil$_2$-groups.
Then one has the following exact sequence:
$$
0\to \qu(G,[H,H])\to \qw(G,H)\to \Hom(G,H\ab).
$$
If additionally $G$ is free in $\Nil$, then the last map is surjective.
\end{Lemma}

\begin{proof} Assume $f:G\to H$ is a q-map. Then the composite of $f$ with
the quotient map $H\to H\ab$ is a homomorphism, which is zero provided
the image of $f$ lies in $[H,H]$. Then the resulting map is quadratic.
Conversely any quadratic map $G\to [H,H]$ considered as a map $G\to H$ is a
q-map. If $G$ is free then any homomorphism
$G\to H\ab$ has a lifting to a homomorphism $G\to H$ and the result follows.
\end{proof}

For any $G\in \Nil$ and any $A\in \Ab$ define the groups $\qu^*(G,A)$ as
the simplicial derived functors of the functor $\qu(-,A)$. More precisely, let
$G_*$ be a free simplicial resolution of $G$. Thus $G_*$ is a simplicial object
in $\Nil$ such that for each $n\ge0$ the group $G_n$ is free in $\Nil$ and
$\pi_i(G_*)=0$ for $i>0$ and $\pi_0(G_*)=G$. Then one can consider the
cosimplicial abelian group $\qu(G_*,A)$. It is well known that the groups
$\pi^*(G_*,A)$ do not depend on the choice of a free simplicial resolution
and they are denoted by $\qu^*(G,A)$. Actually $\qu^0(G,A)=\qu(G,A)$.

For $N\in \Nil$ the sets $\pi^i(\qw(G_*,N))$ for $i=0,1$ also do not depend on
the choice of a free simplicial resolution of $G$. We will denote them by
$\qw^i(G,H)$, $i=0,1$. Actually $\qw^0(G,H)=\qw(G,H)$.

\begin{Proposition}
Let $G$ and $H$ be nil$_2$-groups. Then one has the following exact sequence:
\begin{multline*}
0\to \qu(G,[H,H])\to \qw(G,H)\to \Hom(G,H\ab)\\
\to \qu^1(G,[H,H])\to \qw^1(G,H)\to H^2_\Nil(G,H\ab),
\end{multline*}
where all terms are groups except for $\qw^1(G,H)$ and all maps are homomorphisms
except for the last two maps.
\end{Proposition}

\begin{proof}
By Lemma \ref{quqwzustia} we have a short exact sequence of cosimplicial
groups $$0\to \qu(G_*,[H,H])\to \qw(G_*,H)\to \Hom(G_*,H\ab)\to 0,$$ where
$G_*$ is a free simplicial resolution of $G$. The rest follows from Lemma
\ref{exactnonabelianhomotopy}.
\end{proof}

\begin{Corollary}
For nil$_2$-groups $G$ and $H$ and a homomorphism $f:G\to H\ab$, there is a
well-defined element $o(f)\in\qu^1(G,[H,H])$ which vanishes if and only if $f$
lifts to a q-map $G\to H$.
\end{Corollary}\qed

In particular, taking above $H$ to be arbitrary, $G=H\ab$ and $f$ the identity
map, denote the corresponding element $o(f)$ by $o(H)$; this is thus an
element in $\qu^1(H\ab,[H,H])$. Then we have

\begin{Corollary}
For a nil$_2$-group $G$ there is a
well-defined element $o(G)\in\qu^1(G\ab,[G,G])$ which vanishes if and only if
$G$ is q-split.
\end{Corollary}\qed

\section*{Acknowledgements}

The paper was written during visits of the authors to the University of
Bielefeld and the Max-Planck-Institut f\"ur Mathematik in Bonn. The authors
gratefully acknowledge hospitality of these institutions.

\end{document}